\documentclass[11pt]{amsart}
\usepackage[hyperfootnotes=false]{hyperref}
\usepackage{enumerate}
\usepackage{tikz}
\usepackage{tikz-cd}
\usepackage{hyperref}
\usepackage{amsthm,amsfonts,amsmath,amscd,amssymb}
\usepackage{xcolor,float}
\usepackage[all]{xy}
\usepackage{mathrsfs}
\usepackage{graphics,graphicx}
\usepackage{caption}

\let\oldmarginpar\marginpar
\renewcommand\marginpar[1]{\-\oldmarginpar[\raggedleft\footnotesize #1]%
	{\raggedright\footnotesize #1}}
\theoremstyle{plain}
\newtheorem{thm}{Theorem}[section]
\newtheorem{lemma}[thm]{Lemma}

\newtheorem{prop}[thm]{Proposition}
\newtheorem{cor}[thm]{Corollary}

\theoremstyle{definition}

\newtheorem{definition}[thm]{Definition}
\newtheorem{remark}[thm]{Remark}

\newtheorem{ex}[thm]{Example}
\theoremstyle{remark}

\numberwithin{equation}{section}

\renewcommand{\L}{\mathbb{L}}

\newcommand{\N}{\mathbb{N}}
\newcommand{\Z}{\mathbb{Z}}

\newcommand{\R}{\mathbb{R}}
\newcommand{\C}{\mathbb{C}}

\newcommand{\SH}{\mathcal{H}}

\newcommand{\SM}{\mathcal{M}}

\newcommand{\SF}{\mathscr{F}}

\newcommand{\La}{\Lambda}
\newcommand{\la}{\lambda}

\renewcommand{\a}{\alpha}

\newcommand{\sse}{\subseteq}
\newcommand{\lr}{\longrightarrow}

\newcommand{\Gr}{\operatorname{Gr}}

\newcommand{\Aut}{\operatorname{Aut}}

\newcommand{\Cont}{\operatorname{Cont}}
\newcommand{\SL}{\operatorname{SL}}
\newcommand{\GL}{\operatorname{GL}}
\newcommand{\PSL}{\operatorname{PSL}}

\newcommand{\Sh}{\operatorname{Sh}}
\newcommand{\Aug}{\operatorname{Aug}}

\newcommand{\st}{\text{st}}

\newcounter{daggerfootnote}

\newcommand{\cL}{\mathcal{L}}
\newcommand{\cM}{\mathcal{M}}

\newcommand{\cO}{\mathcal{O}}

\begin{document}

\title{Infinitely many Lagrangian fillings}
\subjclass[2010]{Primary: 53D10. Secondary: 53D15, 57R17.}

\author{Roger Casals}
\address{University of California Davis, Dept. of Mathematics, Shields Avenue, Davis, CA 95616, USA}
\email{casals@math.ucdavis.edu}

\author{Honghao Gao}
\address{Department of Mathematics, Michigan State University, 619 Red Cedar Road, East Lansing, MI 48824, USA}
\email{gaohongh@msu.edu}

\maketitle

\begin{abstract}
We prove that all maximal-tb positive Legendrian torus links $(n,m)$ in the standard contact 3-sphere, except for $(2,m)$, $(3,3),(3,4)$ and $(3,5)$, admit infinitely many Lagrangian fillings in the standard symplectic 4-ball. This is proven by constructing infinite order Lagrangian concordances which induce faithful actions of the modular group $\PSL(2,\Z)$ and the mapping class group $M_{0,4}$ into the coordinate rings of algebraic varieties associated to Legendrian links. In particular, our results imply that there exist Lagrangian concordance monoids with subgroups of exponential-growth, and yield Stein surfaces homotopic to a 2-sphere with infinitely many distinct exact Lagrangian surfaces of higher-genus. We also show that there exist infinitely many satellite and hyperbolic knots with Legendrian representatives admitting infinitely many exact Lagrangian fillings.
\end{abstract}

\section{Introduction}

We show that essentially all maximal-tb positive Legendrian torus links in the standard contact 3-sphere remarkably admit infinitely many non-Hamiltonian isotopic exact Lagrangian fillings in the standard symplectic $4$-ball. Heretofore, the existence of Legendrian links with infinitely many exact Lagrangian fillings remained open.

In fact, the faithful $\PSL(2,\Z)$ representation in our Theorem \ref{thm:main} allows us to obtain several consequences. We present new results for Lagrangian concordance monoids, including the first known example of a Lagrangian concordance of infinite order, the existence of an exponential-growth subgroup in the fundamental group of the space of Legendrian links isotopic to $\La(3,6)$, and the existence of Weinstein 4-manifolds homotopic to the 2-sphere with infinitely many non-Hamiltonian isotopic exact Lagrangian surfaces of higher-genus in the same smooth isotopy class. In addition, we construct infinitely many instances of both satellite and hyperbolic knots in the 3-sphere with Legendrian representatives with infinitely many exact Lagrangian fillings in the standard symplectic 4-ball.

\subsection{Context} Legendrian knots in contact 3--manifolds are instrumental to study the contact geometry of 3--manifolds \cite{Bennequin83,EliashbergFraser98,Etnyre03,Etnyre05,EtnyreHonda01a,EtnyreNgVertesi13,Geiges08}. The classification of Legendrian knots and their Lagrangian fillings has been one of the central areas of research in low-dimensional contact topology \cite{LF6,EkholmHondaKalman16,LF3,PanCatalanFillings,LF4,LF2,LF1}. The only Legendrian knot for which there exists a complete non-empty classification of Lagrangian fillings is the Legendrian unknot \cite{EliashbergPolterovich96}.

The works \cite{EkholmHondaKalman16,PanCatalanFillings,STWZ} succeeded in constructing a Catalan number worth of Lagrangian fillings for the maximal-tb positive Legendrian $(2,n)$-torus links. It is also known that all positive braids admit at least one Lagrangian filling \cite{Kalman-braid}; see also \cite{EkholmHondaKalman16,LF3}. A crucial question that remained open is the existence of Legendrian links with infinitely many exact Lagrangian fillings. This article affirmatively resolves this question.

\begin{center}
	\begin{figure}[h!]
		\centering
		\includegraphics[scale=0.7]{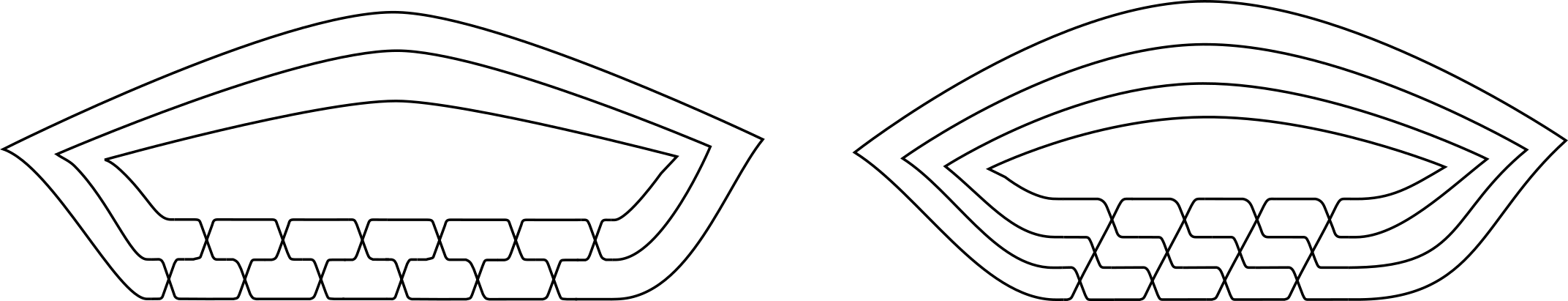}
		\caption{The Legendrian torus links $\La(3,6)$ (left) and $\La(4,4)$ (right).}
		\label{fig:BraidsIntro}
	\end{figure}
\end{center}
In fact, we shall geometrically construct Lagrangian concordances which themselves produce infinitely many Lagrangian fillings, which is a significantly stronger statement than the existence of infinitely many Lagrangian fillings. The constructions are explicit and can be readily drawn in the front projection. The construction implies that these exact Lagrangian surfaces are all smoothly isotopic. We will distinguish these Lagrangian fillings by studying their action on part of the coordinate ring of the moduli of framed constructible sheaves $\SM(\La)$ \cite{GKS_Quantization,KashiwaraSchapira_Book,STZ} for certain Legendrian links $\La\sse(S^3,\xi_\st)$. The techniques we use for our results illustrate the strength of applying methods from the microlocal theory of sheaves \cite{KashiwaraSchapira_Book,STZ} and the theory of cluster algebras \cite{FockGoncharov_ModuliLocSys,FominZelevinsky_ClusterI,Fraser} to 3-dimensional contact and symplectic topology.

\subsection{Main Results} Let $\La(n,m)\sse(S^3,\xi_\st)$ be the maximal-tb positive Legendrian $(n,m)$-torus link, $(n,m)\in\N\times\N$, as depicted in Figure \ref{fig:BraidsIntro}. Positive Legendrian torus knots are Legendrian simple positive braids \cite{EtnyreHonda01b}, and thus are uniquely determined by their Thurston-Bennequin invariants and their rotation numbers. These Legendrian links can be obtained by considering the positive braids $\beta=(\sigma_1\sigma_2\cdot\ldots\cdot\sigma_{n-1})^{m+n}$ in $(J^1S^1,\xi_\st)$ and satelliting the zero section $S^1\sse(J^1S^1,\xi_\st)$ to the standard Legendrian unknot $\La(1,1)\sse(S^3,\xi_\st)$. Let $\cL(n,m)$ be the space of Legendrian links isotopic to the maximal-tb Legendrian torus link $\La(n,m)\sse(S^3,\xi_\st)$, with base point an arbitrary but fixed maximal-tb Legendrian representative.

Let $\cM(\La(n,m))$ be the moduli space of framed sheaves associated to the Legendrian link $\La(n,m)$, as we shall introduce in Section \ref{sec:algebra}. This $\cM(\La(n,m))$ is an algebraic variety \cite{STWZ}, and in our case it will be a quasi-projective subvariety of the projective Grassmannian $\Gr(n,n+m)$. Since $\cM(\La(n,m))$ is a Legendrian isotopy invariant \cite{GKS_Quantization,STWZ,STZ}, it defines a monodromy representation
$$\Gamma:\pi_1(\cL(n,m))\lr\Aut(\cM(\La(n,m))),$$
into the space of algebraic automorphisms of the algebraic variety $\cM(\La(n,m))$. In turn, by pull-back, we obtain a representation
$$\Gamma^*:\pi_1(\cL(n,m))\lr\Aut(\C[\cM(\La(n,m))]),$$
into the automorphisms of the coordinate ring $\C[\cM(\La(n,m))]$ of $\cM(\La(n,m))$. In particular, a set of based loops $C_1,\ldots,C_r:S^1\lr\cL(n,m)$, $r\in\N$, gives rise to a monodromy representation
$$\Gamma^*:\langle [C_1],\ldots,[C_r]\rangle \lr\Aut(\C[\cM(\La(n,m))]).$$
of the subgroup $\langle [C_1],\ldots,[C_r]\rangle\leq\pi_1(\cL(n,m))$ generated by the homotopy classes of the based maps $C_1,\ldots,C_r:S^1\lr\cL(n,m)$. The first result we present is:

\begin{thm}\label{thm:main}
Let $\cL(3,6)$ be the space of Legendrian links isotopic to the maximal-tb Legendrian torus link $\La(3,6)\sse(S^3,\xi_\st)$. Then there exist two based loops $A,B:S^1\longrightarrow\cL(3,6)$, and a regular function $\Delta\in\C[\cM(\La(3,6))]$ such that the monodromy representation
$$\Gamma^*:\langle [A],[B] \rangle\longrightarrow \Aut(\C[\cM(\La(3,6))])$$
restricts to a faithful modular representation
$$\Gamma^*|_{\cO(\Delta)}:\PSL(2,\Z) \longrightarrow \Aut(\cO(\Delta))$$
along the orbit $\cO(\Delta)$ of the function $\Delta$.\hfill$\Box$
\end{thm}

In Theorem \ref{thm:main}, we choose the base point for the space $\cL(3,6)$ to be the Legendrian link in $(\R^3,\xi_\st)$ whose front projection is depicted in Figure \ref{fig:BraidsIntro}, under an arbitrary but fixed choice of Darboux chart $(\R^3,\xi_\st)\sse(S^3,\xi_\st)$. In the statement of Theorem \ref{thm:main}, $\langle [A],[B] \rangle\sse\pi_1(\cL(3,6))$ denotes the subgroup generated by the homotopy classes $[A],[B]\in\pi_1(\cL(3,6))$, with concatenation of based loops as its group operation. The modular group $\PSL(2,\Z)$ shall appear geometrically as the free product $\Z_3*\Z_2$, with the factor $\Z_3$ generated by the restriction of $[A]$ and the factor $\Z_2$ generated by the restriction of $[B]$, and $\Aut(\C[\cM(\La(3,6))])$ denotes the group of (cluster) automorphisms of $\C[\cM(\La(3,6))]$. Theorem \ref{thm:main} is remarkable in that $\PSL(2,\Z)$ is an infinite group and thus provides the first result of its kind in the study of $3$-dimensional Legendrian links.

\begin{remark}\label{rmk:J10} The reason for the choice of the Legendrian torus link $\La(3,6)$ is that it is the geometric source of the extended root system $E_8^{(1,1)}$. Indeed, it can be understood as the maximal-tb Legendrian approximation of the transverse link of the unimodal $J_{10}$ singularity \cite{Arnold90}. The proof of Theorem \ref{thm:main} shall clarify how the $E_8^{(1,1)}$ algebraic structure arises from $\La(3,6)$.\hfill$\Box$
\end{remark}

We show that the Legendrian torus link $\La(4,4)$ also exhibits a noteworthy symmetry:

\begin{thm}\label{thm:main2}
	Let $\cL(4,4)$ be the space of Legendrian links isotopic to the maximal-tb Legendrian torus link $\La(4,4)\sse(S^3,\xi_\st)$. Then there exist three based loops $\Xi_i:S^1\longrightarrow\cL(4,4)$, $1\leq i\leq 3$, and a subset $F\sse\C[\cM(\La(4,4))]$ such that the monodromy representation
	$$\Gamma^*:\langle [\Xi_1],[\Xi_2],[\Xi_3] \rangle\longrightarrow \Aut(\C[\cM(\La(4,4))])$$
restricts to a faithful representation
$$\Gamma^*|_{\cO(F)}:M_{0,4} \longrightarrow \Aut(\cO(F))$$
of the mapping class group $M_{0,4}$ along the orbit $\cO(F)$.\hfill$\Box$
\end{thm}

The mapping class group $M_{0,4}$ of the four-punctured 2-sphere contains a subgroup isomorphic to $\PSL(2,\Z)$ with finite index and it is thus infinite. In Theorem \ref{thm:main2}, the base point for the space $\cL(4,4)$ is chosen to be the Legendrian link in $(\R^3,\xi_\st)$ with front projection as depicted in Figure \ref{fig:BraidsIntro}, also under an arbitrary but fixed choice of Darboux chart $(\R^3,\xi_\st)\sse(S^3,\xi_\st)$. The subgroup $\langle [\Xi_1],[\Xi_2],[\Xi_3] \rangle\sse\pi_1(\cL(4,4))$ has loop concatenation as its group operation.

\begin{remark}
The two groups $\PSL(2,\Z)$ and $M_{0,4}$ featured in Theorems \ref{thm:main} and \ref{thm:main2} are akin to each other in that there are two group isomorhisms $\PSL(2,\Z)\cong B_3/Z(B_3)$ and $M_{0,4}\cong B^s_4/Z(B^s_4)$, where $B_3$ denotes the braid group in 3-strands, $B^{s}_4$ denotes the spherical braid group in 4-strands, and $Z(B_3)$ and $Z(B^s_4)$ denote their respective centers.\hfill$\Box$
\end{remark}

Let us now state implications of Theorems \ref{thm:main} and \ref{thm:main2}, all of which are new results in low-dimensional contact and symplectic topology.

\subsection{Lagrangian Fillings} Consider the subset
$$\SH:=\{(n,m)\in\N\times\N:n\leq m, 3\leq n,6\leq m\}\cup\{(4,4),(4,5),(5,5)\}\sse \N\times\N.$$
The first consequence of Theorems \ref{thm:main} and \ref{thm:main2} is:

\begin{cor}\label{cor:infinitefillings}
The Legendrian torus link $\La(n,m)\sse (S^3,\xi_\st)$, $(n,m)\in\SH$, admits infinitely many exact Lagrangian fillings in the standard symplectic $4$-ball $(D^4,\omega_\st)$.\hfill$\Box$
\end{cor}

For each $\La(n,m)$, these infinitely many exact Lagrangian fillings are smoothly isotopic and not Hamiltonian isotopic. Note that both Theorem \ref{thm:main} and Theorem \ref{thm:main2} are needed in order to cover all $\La(n,m)$ for $(n,m)\in\SH$. That said, Theorem \ref{thm:main} suffices in order to conclude Corollary \ref{cor:infinitefillings} for $(n,m)\in(\SH\setminus \{(4,4),(4,5),(5,5)\})$ and thus Theorem \ref{thm:main2} is included to achieve Corollary \ref{cor:infinitefillings} for $(n,m)=(4,4),(4,5)$ and $(5,5)$. It should be noted that the article \cite{STWZ} succeeded in constructing {\it finitely} many Lagrangian fillings of the maximal-tb Legendrian $(n,m)$-torus link, as many as maximal pairwise weakly separated $n$-element subsets \cite{OPS,Scott} of $[1,n+m]$, a finite number which is bounded above by $nm+1$. Corollary \ref{cor:infinitefillings} implies that these finitely many exact Lagrangian fillings do not exhaust all possible, actually infinitely many, exact Lagrangian fillings.

Every knot $K\sse S^3$ is either a torus knot, a satellite knot or a hyperbolic knot, as proven in \cite[Theorem 2.3]{Thurston_Knots} by W.P.~Thurston. Let us consider a Legendrian representative $\La_K\sse(S^3,\xi_\st)$ of the smooth type $K\sse S^3$ and denote by $l(\La_K)\in\N\cup\{\infty\}$ the number of orientable exact Lagrangian fillings $L\sse (D^4,\omega_\st)$ of the Legendrian knot $\La_K\sse(S^3,\xi_\st)$, up to a Hamiltonian isotopy. Consider the smooth invariant
$$\overline{l}(K):=\sup\{l(\La_K):\La_K\sse(S^3,\xi_\st)\mbox{ is a Legendrian representative of } K\}\in\N\cup\{\infty\}$$
for a smooth knot $K\sse S^3$. To our knowledge, there is no hitherto known instance of a non-trivial knot $K\sse S^3$ for which $\overline{l}(K)$ is known and non-vanishing. In addition, there are non-trivial knots $K\sse S^3$ for which the invariant $\overline{l}(K)$ vanishes. For instance, $\overline{l}(m(8_{19}))=0$ is known to vanish since the Kauffman upper bound is not sharp \cite{FuchsTabachnikov_Kauffman,Rutherford_Kauffman}. We shall now use Theorem \ref{thm:main} to show that $\overline{l}(K)$ is actually infinite for infinitely many knots within each of the three Thurston classes:

\begin{cor}\label{cor:lisinfinite} The equality $\overline{l}(K)=\infty$ holds for infinitely many torus knots, infinitely many satellite knots and infinitely many hyperbolic knots $K\sse S^3$ in the $3$-sphere.
\end{cor}

The satellite knots in Corollary \ref{cor:lisinfinite} can be chosen to be cable knots, and the hyperbolic knots we will exhibit are also well-beloved \cite{Birman_Lorenz,Ghys_Knots}. For instance, we will show that $\overline{l}(K)=\infty$ already for $K=k(4_3)$, one of the simplest hyperbolic non-2-bridge knots \cite{CDW_SimplestHyperbolic}.

\subsection{Lagrangian Concordances} Now, let $\L(n,m)$ be the monoid of exact Lagrangian concordances, up to Hamiltonian isotopy, for the Legendrian link $\La(n,m)\sse (S^3,\xi_\st)$. Theorems \ref{thm:main} and \ref{thm:main2} readily imply:

\begin{cor}\label{cor:concordances}
There exists subgroups $\Gamma\sse\L(3,6)$ and $\Gamma'\sse\L(4,4)$ such that the group $\PSL(2,\Z)$ is a factor of $\Gamma$, and $M_{0,4}$ is a factor of $\Gamma'$.\hfill$\Box$
\end{cor}

By definition, the groups $\Gamma$ and $\Gamma'$ in Corollary \ref{cor:concordances} are the subgroups generated by the exact Lagrangian concordances obtained by graphing the Legendrian loops in Theorems \ref{thm:main} and \ref{thm:main2}. Corollary \ref{cor:concordances} emphasizes the relevance of Lagrangian concordances in the study of Legendrian knots. In particular, the existence of Lagrangian concordances of infinite order is a new result which itself provides a genuinely useful perspective for the study of Lagrangian fillings. Indeed, there is no a priori reason for the infinite Lagrangian fillings in Corollary \ref{cor:infinitefillings} to be describable in terms of a finite number of Lagrangian concordances. The present results show that this is the case. Similarly,

\begin{cor}\label{cor:fundamentalgroup}
	There exists subgroups $\Gamma\sse\pi_1(\cL(3,6))$ and $\Gamma'\sse\pi_1(\cL(4,4))$ such that the group $\PSL(2,\Z)$ is a factor of $\Gamma$, and $M_{0,4}$ is a factor of $\Gamma'$.\hfill$\Box$
\end{cor}

Corollary \ref{cor:concordances} and Corollary \ref{cor:fundamentalgroup} are the first instances in contact topology of infinite order elements in the concordance monoid $\L(\La)$, and the fundamental group $\pi_1(\cL(\La))$, for a Legendrian $\La\sse(S^3,\xi_\st)$. Both $\PSL(2,\Z)$ and $M_{0,4}$ contain free groups of any countable rank as subgroups, and thus many infinite order elements exist in $\pi_1(\cL(3,6))$ and $\pi_1(\cL(4,4))$. In fact, $\Gamma$ and $\Gamma'$ are exponential-growth subgroups of $\pi_1(\cL(3,6))$ and $\pi_1(\cL(4,4))$.

\begin{remark}
Corollary \ref{cor:fundamentalgroup} stands in contrast with A. Hatcher's work \cite[Theorem 1]{Hatcher_Knots} in the smooth category. Indeed, the fundamental group $\pi_1(\mathscr{K}_{(n,m)})$ of the space $\mathscr{K}_{(n,m)}$ of smooth knots in $S^3$ isotopic to the $(n,m)$-torus knot is the finite Abelian group $\Z_2$.\hfill$\Box$
\end{remark}

\subsection{Stein surfaces} Finally, let $M(n,m)$ be the Stein surface obtained by a attaching Weinstein $2$-handle \cite{CieliebakEliashberg12,Weinstein91} to $(D^4,\la_\st)$ along each of the components of the Legendrian link $\La(n,m)\sse(S^3,\xi_\st)$. For $\gcd(n,m)=1$, the Weinstein 4-manifold $M(n,m)$ is homotopy equivalent to the 2-sphere. Theorems \ref{thm:main} and \ref{thm:main2} imply the existence of infinitely many Lagrangian surfaces in the following Stein surfaces:
\begin{cor}\label{cor:weinstein}
Let $(n,m)\in\SH$, $\gcd(n,m)=1$ and $M(n,m)$ the Weinstein 4-manifold obtained by attaching a Weinstein 2-handle to $(D^4,\omega_\st)$ along $\La(n,m)$. Then $M(n,m)$ contains infinitely many smoothly isotopic closed exact Lagrangian surfaces of genus $\frac{1}{2}(n-1)(m-1)$ which are not Hamiltonian isotopic.\hfill$\Box$
\end{cor}
To our knowledge, Corollary \ref{cor:weinstein} presents the first known Stein surfaces homotopic to the $2$-sphere with infinitely many non-Hamiltonian isotopic exact Lagrangian surfaces of higher genus in the same smooth isotopy class.

Infinitely many distinct Lagrangian 2-spheres were known to exist in $A_k$-Milnor fibres \cite[Theorem 5.10]{Seidel_Graded}, $k\geq3$, and infinitely many exact Lagrangian tori were known to exist in certain Stein surfaces by using either of the articles \cite{Keating,STW_Combinatorics,Vianna2,Vianna1}. (These infinite families of genus $0$ and $1$ are presently not known to come from infinitely many Lagrangian fillings of a Legendrian link, nor are the ambient Weinstein 4-manifolds homotopic to the $2$-sphere.)

In Corollary \ref{cor:weinstein}, the 1-dimensional intersection form of the Weinstein 4-manifold $M(n,m)$ is positive definite and equals $(nm-n-m-1)$, since $tb(\La(n,m))=nm-n-m$. In consequence, $M(n,m)$ does not admit any Lagrangian surface of genus strictly less than $\frac{1}{2}(n-1)(m-1)$. Thus, the genus in Corollary \ref{cor:weinstein} is sharp.

\begin{remark}\label{rmk:Lagrsph}
The Lagrangian 2-spheres in \cite{Seidel_Graded} differ by a composition of symplectic Dehn twists \cite{Arnold95,Seidel99}. This is not the case for the exact Lagrangian higher-genus surfaces in Corollary \ref{cor:weinstein} since, by the paragraph above, the Weinstein 4-manifolds $M(n,m)$, $(n,m)\in\SH$, do not contain embedded Lagrangian 2-spheres.\hfill$\Box$
\end{remark}
{\bf Organization.} The article is organized as follows. Section \ref{sec:loops} geometrically constructs the loops in Theorem \ref{thm:main} and Theorem \ref{thm:main2}. Section \ref{sec:algebra} provides the necessary aspects from the theory of Legendrian invariants constructed through the study of microlocal sheaves. Sections \ref{sec:mainproof1} and \ref{sec:mainproof2} prove Theorem \ref{thm:main} and Theorem \ref{thm:main2}, respectively, and Section \ref{sec:corollaries} proves the Corollaries stated in the introduction.\hfill$\Box$\\

{\bf Acknowledgements.} R.C. is grateful to J.B. Etnyre and L. Ng for useful conversations, and to I. Smith for helpful comments on Corollary \ref{cor:weinstein}. J.B. Etnyre asked us an interesting question on our first manuscript, now answered in Corollary \ref{cor:lisinfinite}. We also thank A. Keating, J. Sabloff, L. Starkston and U. Varolg\"une\c{s} for their interest on this project and a number of useful comments and suggestions. Finally, both authors are indebted to E. Zaslow for many valuable discussions on microlocal Legendrian invariants, and to C. Fraser for helpful conversations on cluster modular groups. We also thank the referee for their thorough comments and suggestions. R.~Casals is supported by the NSF CAREER grant DMS-1942363 and a Sloan Research Fellowship of the Alfred P. Sloan Foundation.\hfill$\Box$

%
%

\section{The geometric construction}\label{sec:loops}

In this section we construct the Legendrian loops in Theorems \ref{thm:main} and \ref{thm:main2} associated to the Legendrian links $\La(3,6)$ and $\La(4,8)$. This construction is one of the central geometric contributions of the article. This section also serves to setup the elements of contact geometry that we shall need \cite{Etnyre03,Geiges08}.

The Legendrian loops $\Sigma_1,\delta^2,\Xi_1,\Xi_2$ and $\Xi_3$ that we construct can be equivalently considered as exact Lagrangian concordances in the symplectization $(S^3\times\R(t),d(e^t\a_\st))$ with no critical points with respect to the projection onto the $\R$-factor \cite{CieliebakEliashberg12,Geiges08}. These Lagrangian concordances are obtained by graphing concatenations of the Legendrian isotopies described in Subsection \ref{ssec:Legloop}.

\subsection{The standard contact 3-space} The Legendrian links $\La\sse(S^3,\xi_\st)$ in this article will be considered inside the standard contact 3-space $(\R^3,\xi_\st)$, considered as a standard Darboux chart within the contact 3-sphere $(S^3,\xi_\st)$ \cite{Arnold90,ArnoldGivental01}.

In discussing Lagrangian fillings, the inclusion $\La\sse(\R^3,\xi_\st)$ will be composed with the inclusion $(\R^3,\xi_\st)\sse(S^3,\xi_\st)$ given by the one-point compactification. In this identification, the Lagrangian fillings of a Legendrian link $\La\sse(\R^3,\xi_\st)$ will be exact Lagrangian surfaces in $(D^4,\la_\st)$ considered up to Hamiltonian isotopy.

The constructions in this article give rise to contact geometric objects in $(\R^3,\xi_\st)$, including Legendrian links and contact isotopies. Nevertheless, it is enlightening to focus on a small neighborhood of the standard Legendrian unknot $\La_{un}\sse(\R^3,\xi_\st)$ which is contactomorphic to $(J^1S^1,\xi_\st)$, and work in the solid torus $(J^1S^1,\xi_\st)$. Thus, in this article, Legendrian links $\La\sse(J^1S^1,\xi_\st)$ and compactly supported contact isotopies in $\{\varphi_t\}_{t\in[0,1]}\in\Cont(J^1S^1,\xi_\st)$ shall implicitly be understood as Legendrian links $\La\sse(\R^3,\xi_\st)$ and contact isotopies $\{\varphi_t\}_{t\in[0,1]}\in\Cont(\R^3,\xi_\st)$ by satelliting the zero section $S^1\sse(J^1S^1,\xi_\st)$ to $\La_{un}\sse(\R^3,\xi_\st)$.

\subsection{Legendrian Loops}\label{ssec:Legloop} By definition, a Legendrian loop in $(J^1S^1,\xi_\st)$ is a Legendrian isotopy $\{\La_t\}_{t\in[0,1]}\sse(J^1S^1,\xi_\st)$ such that $\La_0=\La_1$. Let $(\theta,p_\theta,z)\in S^1\times\R^2$ be global coordinates in $J^1S^1$ and $\xi_\st=\ker{(dz-p_{\theta}d\theta)}$. The description of our Legendrian loops shall use the front projection
$$(J^1S^1,\xi_\st)\lr S^1\times\R,\quad (\theta,p_\theta,z)\longmapsto(\theta,z),$$
which is indeed a valid front as the fibers are Legendrians. By definition, the Legendrian $\La(\mathscr{B})\sse(\R^3,\xi_\st)$ associated to a positive Legendrian braid $\mathscr{B}\sse(J^1S^1,\xi_\st)$ is the image of $\mathscr{B}$ under the operation of satelliting the zero section $S^1\sse(J^1S^1,\xi_\st)$ along the standard Legendrian unknot.

Let $k\in\N$, a geometric positive braid $\mathscr{B}\sse (J^1S^1,\xi_\st)$, which is a Legendrian link, can be encoded algebraically by a positive expression $\beta$,  i.e. a positive braid word, of an element $[\beta]\in B_k$ of the $k$-stranded braid group
$$B_k:=\langle \sigma_1,\ldots,\sigma_{k-1}|\sigma_i\sigma_{i+1}\sigma_i=\sigma_{i+1}\sigma_i\sigma_{i+1},\sigma_i\sigma_j=\sigma_j\sigma_i\mbox{ for }j\neq i\pm 1,1\leq i,j\leq k-1\rangle,$$
where $\sigma_i$ are the standard Artin generators. The choice of representation $\beta$, i.e. braid word, for the element $[\beta]\in B_k$ is not unique, as one might use the word relations in $B_k$ to obtain different representations $\beta_1,\beta_2$ such that $[\beta_1]=[\beta_2]$ in $B_k$. Given a positive braid word $\beta$, we denote by $\mathscr{B}(\beta)\sse (J^1S^1,\xi_\st)$ the Legendrian link associated to $\beta$, and denote by $\La(\beta)\sse(\R^3,\xi_\st)$ the Legendrian link $\La(\mathscr{B}(\beta))$.

Notice that the geometric braid $\mathscr{B}\sse (J^1S^1,\xi_\st)$ has a front in $S^1\times\R$. Thus, we fix a basepoint $\theta_0\in S^1$ and require that a braid word $\beta$ representing $\mathscr{B}$ has the form:
$$\beta=\prod_{j=1}^{l}\sigma_{i_j},\quad 1\leq i_j\leq k-1,$$
where $\sigma_{i_1}$ is the first crossing in the front diagram of $\mathscr{B}$ on the right of the vertical line $\{\theta_0\}\times\R\sse S^1\times\R$ and the crossings are read from left to right.

In this article, we shall construct Legendrian loops by performing Legendrian isotopies which primarily consist of Reidemeister moves in the front. In particular, the three central operations that we use are:

\begin{itemize}
	\item[(i)] {\bf Reidemeister III Moves}. In terms of the given braid word presentation $\mathscr{B}=\mathscr{B}(\beta)$, the Reidemeister III move consists in applying the relation $\sigma_i\sigma_{i+1}\sigma_i=\sigma_{i+1}\sigma_i\sigma_{i+1}$. We shall refer to a Legendrian isotopy which implements the substitution
	
$$\sigma_i\sigma_{i+1}\sigma_i\longmapsto\sigma_{i+1}\sigma_i\sigma_{i+1}$$

\noindent as an ascending Reidemeister III Move, and denote it by R3$^a$. Similarly, to a Legendrian isotopy which implements the substitution

$$\sigma_{i+1}\sigma_i\sigma_{i+1}\longmapsto\sigma_i\sigma_{i+1}\sigma_i$$

\noindent as a descending Reidemeister III Move, and denote it by R3$^d$. Thus, either R3 Reidemeister move is understood as a Legendrian isotopy.\\

	\item[(ii)] {\bf Cyclic Permutation}. Consider a braid $\mathscr{B}(\beta)$ represented by
	$$\beta=\prod_{j=1}^{l}\sigma_{i_j},\quad 1\leq i_j\leq k-1.$$
	By definition, a cyclic shift $\delta$ is a Legendrian isotopy $\{\psi_t\}_{t\in[0,1]}$ which brings the geometric braid $\mathscr{B}(\beta)$ to $\psi_1(\mathscr{B}(\beta))$ such that the braid word $\psi_1(\beta)$ for the latter is
	$$\psi_1(\beta)=\left(\prod_{j=2}^{l}\sigma_{i_j}\right)\sigma_{i_1},\quad 1\leq i_j\leq k-1.$$
	Note that this braid word for $\psi_1(\beta)$ is read with respect to the fixed basepoint $\theta_0$. Explicitly, this Legendrian isotopy can be geometrically visualized by rotating $\mathscr{B}(\beta)$ to the left by an appropriate angle while keeping the zero section fixed.\\
	
	\noindent Since we study Legendrian braids in $(J^1S^1,\xi_\st)$, rather than $(J^1[0,1],\xi_\st)$, two braid words $\beta_1=\beta_2$ which differ by a cyclic permutation yield Legendrian isotopic $\La(\beta_1)$ and $\La(\beta_2)$. Hence, the operations above produce Legendrian isotopies.\\
	
	\item[(iii)] {\bf Commutation}. The third move $\gamma$ is just implementing the commutation relation in the braid group $B_k$. It is described as follows:
	$$\left(\prod_{j=1}^{p-1}\sigma_{i_j}\right)\sigma_{i_p}\sigma_{i_{p+1}}\left(\prod_{j=p+2}^{l}\sigma_{i_j}\right)\longmapsto\left(\prod_{j=1}^{p-1}\sigma_{i_j}\right)\sigma_{i_{p+1}}\sigma_{i_p}\left(\prod_{j=p+2}^{l}\sigma_{i_j}\right),$$
	with indices $1\leq i_j\leq k-1, 1\leq p\leq l-1$ and $i_{p+1}\neq i_{p}\pm1$.  This move can be realized by a compactly supported Legendrian isotopy in $(J^1S^1,\xi_\st)$ which we also refer to as $\gamma$, which is the greek letter for $c$, standing for commutation.
	
\end{itemize}

\begin{ex}\label{ex:Kalman} Consider the braid word $\beta=(\sigma_1\sigma_2\cdot\ldots\cdot\sigma_{n-1})^{m+n}$ which geometrically represents the Legendrian torus link $\La(n,m)=\La(\beta)\sse(\R^3,\xi_\st)$. Then the composition of the cyclic shift $\delta$ exactly $(n-1)$-times yields a Legendrian loop $\delta^{n-1}$ for $\La(n,m)$. This is the Legendrian loop studied in \cite{Kalman} where it is shown to be a non-trivial Legendrian loop. We shall provide our own alternative proof of this non-triviality.\hfill$\Box$
\end{ex}

\subsection{The $\Sigma_1$-Loop for $\La(3,6)$} In this subsection we define a Legendrian loop $\Sigma_1$ for the maximal-tb Legendrian links $\La(3,3s)\sse(\R^3,\xi_\st)$, $s\in\N$, represented by the positive braid $\mathscr{B}(\beta)$, with braid word $\beta=(\sigma_1\sigma_2)^{3(s+1)}$, in the front domain $S^1\times\R$. The loop $\Sigma_1$ is defined as the composition of Legendrian isotopies induced by the following sequence of moves:
\begin{align*}
(\sigma_1\sigma_2)^{3(s+1)}&=(\sigma_1\sigma_2\sigma_1\sigma_2\sigma_1\sigma_2)^{s+1} \stackrel{\delta}{\approx} (\sigma_2\sigma_1\sigma_2\sigma_1\sigma_2\textcolor{blue}{\underline{\sigma_1}})^{s+1}\stackrel{\mbox{\tiny R3}^d}{\approx} (\textcolor{blue}{\underline{\sigma_1\sigma_2\sigma_1}}\sigma_1\sigma_2\sigma_1)^{(s+1)}\\
&\stackrel{\mbox{\tiny R3}^a}{\approx} (\sigma_1\sigma_2\sigma_1\textcolor{blue}{\underline{\sigma_2\sigma_1\sigma_2}})^{(s+1)}=(\sigma_1\sigma_2)^{3(s+1)}.
\end{align*}

In the above sequence, the underlined letters in blue indicate changes in the braid word. In words, the first isotopy is a cyclic shift moving $\sigma_1$ to the end of the braid by shifting left past $\{\theta_0\}\times\R\sse S^1\times\R$. The second isotopy consists of $(3s+3)$ simultaneous and commuting Reidemeister R3$^d$ moves, whereas the third isotopy consists of $(3s+3)$ simultaneous and commuting Reidemeister R3$^a$ moves. The composition of these isotopies yields the initial braid word $(\sigma_1\sigma_2)^{3(s+1)}$ and thus it generates a Legendrian loop. 

\begin{definition}\label{def:Sigma1Loop} Consider $\La(3,3s)\sse(\R^3,\xi_\st)$, the Legendrian isotopy $\Sigma_1$ is the Legendrian loop of $\La(3,3s)$ induced by the sequence of Legendrian isotopies
$$(\sigma_1\sigma_2)^{3(s+1)}\stackrel{\delta}{\approx} \sigma_2(\sigma_1\sigma_2)^{3s+2}\textcolor{blue}{\underline{\sigma_1}}\stackrel{\mbox{\tiny R3}^d}{\approx} ((\textcolor{blue}{\underline{\sigma_1\sigma_2\sigma_1}})(\sigma_1\sigma_2\sigma_1))^{(s+1)}\stackrel{\mbox{\tiny R3}^a}{\approx} ((\sigma_1\sigma_2\sigma_1)(\textcolor{blue}{\underline{\sigma_2\sigma_1\sigma_2}}))^{(s+1)},$$
once the zero section $S^1\sse(J^1S^1,\xi_\st)$ is satellited to the standard unknot.\hfill$\Box$
\end{definition}

Definition \ref{def:Sigma1Loop} yields Legendrian loops for $\La(3,3s)$ for any $s\in\N$. In this article it shall suffice to focus on the case $s=2$. It might be relevant to notice that in Section \ref{sec:mainproof1} we shall prove that the loop $\Sigma_1$ is non-trivial as a Legendrian loop and it is different from the cyclic loop in Example \ref{ex:Kalman}. In fact, the $\Sigma_1$-loop and the cyclic shift $\delta$ for the braid $\beta=(\sigma_1\sigma_2)^9$ will suffice in order to construct the representation in Theorem \ref{thm:main}.

\begin{remark} The Legendrian loop $\Sigma_1$ is geometrically constructed in order to algebraically act as the first Artin generator for a braid group action of $B_3$ into $\C[\SM(\La(3,6))]$.\hfill$\Box$
\end{remark}

\subsection{The $\Xi_1$-Loop for $\La(4,4)$}\label{ssec:Xi1} Let us now define a Legendrian loop $\Xi_1$ for the maximal-tb Legendrian links $\La(4,4s)\sse(\R^3,\xi_\st)$, $s\in\N$, represented by the 4-stranded positive braid $\mathscr{B}(\beta)$, with braid word $\beta=(\sigma_1\sigma_2\sigma_3)^{4(s+1)}$.

The Legendrian loop $\Xi_1$ is described by the cyclic shift
\begin{align}\label{eq:sigma1}(\sigma_1\sigma_2\sigma_3)^{4(s+1)} &\stackrel{\delta}{\approx} \sigma_2\sigma_3(\sigma_1\sigma_2\sigma_3)^{4s+3}\textcolor{blue}{\underline{\sigma_1}}=(\sigma_2\sigma_3\sigma_1)^{4(s+1)}
\end{align}
followed by the sequence of moves:
\begin{align*}
&(\sigma_2\sigma_3\sigma_1\sigma_2\sigma_3\sigma_1\sigma_2\sigma_3\sigma_1\sigma_2\sigma_3\sigma_1)^{(s+1)} \stackrel{\gamma}{\approx} (\sigma_2{\color{blue}{\underline{\sigma_1\sigma_3}}}\sigma_2\sigma_3\sigma_1 \sigma_2\sigma_3\sigma_1\sigma_2\sigma_3\sigma_1)^{(s+1)} \\
&\stackrel{R3^d}{\approx} (\sigma_2\sigma_1{\color{blue}{\underline{\sigma_2\sigma_3\sigma_2}}}\sigma_1 \sigma_2\sigma_3\sigma_1\sigma_2\sigma_3\sigma_1)^{(s+1)} \stackrel{R3^d}{\approx} ({\color{blue}{\underline{\sigma_1\sigma_2\sigma_1}}}\sigma_3\sigma_2\sigma_1 \sigma_2\sigma_3\sigma_1\sigma_2\sigma_3\sigma_1)^{(s+1)}  \quad (\Psi^{(1)}_t \textrm{ until here})\\
&\stackrel{\gamma}{\approx} (\sigma_1\sigma_2{\color{blue}{\underline{\sigma_3\sigma_1}}}\sigma_2\sigma_1 \sigma_2\sigma_3\sigma_1\sigma_2\sigma_3\sigma_1)^{(s+1)} \stackrel{\gamma}{\approx}  (\sigma_1\sigma_2\sigma_3\sigma_1\sigma_2\sigma_1\sigma_2\sigma_3\sigma_1 \sigma_2{\color{blue}{\underline{\sigma_1\sigma_3}}})^{(s+1)} \\
&\stackrel{R3^a}{\approx} (\sigma_1\sigma_2\sigma_3\sigma_1\sigma_2\sigma_1\sigma_2\sigma_3{\color{blue}{\underline{\sigma_2\sigma_1\sigma_2}}}\sigma_3)^{(s+1)} \stackrel{R3^a}{\approx}  (\sigma_1\sigma_2\sigma_3\sigma_1\sigma_2\sigma_1{\color{blue}{\underline{\sigma_3\sigma_2\sigma_3}}}\sigma_1\sigma_2\sigma_3)^{(s+1)} \\
& \stackrel{\gamma}{\approx}  (\sigma_1\sigma_2\sigma_3\sigma_1\sigma_2{\color{blue}{\underline{\sigma_3\sigma_1}}}\sigma_2\sigma_3\sigma_1\sigma_2\sigma_3)^{(s+1)} = (\sigma_1\sigma_2\sigma_3)^{4(s+1)} \quad (\Psi^{(2)}_t \textrm{ until here}).
\end{align*}

In each of the above rows, the underlined letters emphasized in color blue represent those braid generators, equivalently crossings of the front, which have been affected at each step when performing the indicated Legendrian isotopy, consisting either of a Reidemeister R3 move or a cyclic shift $\delta$. Note that the sequence above ends with the braid word $(\sigma_1\sigma_2\sigma_3)^{4(s+1)}$, and thus yields a Legendrian loop when preconcatenated with the Legendrian isotopy in Equation \ref{eq:sigma1}.

\begin{definition}\label{def:Xi1Loop} Consider $\La(4,4s)\sse(\R^3,\xi_\st)$, the Legendrian isotopy $\Xi_1$ is the Legendrian loop of $\La(4,4s)$ given by concatenating the two sequences above and satelliting the zero section $S^1\sse(J^1S^1,\xi_\st)$ to the standard unknot.\hfill$\Box$
\end{definition}

Let us now proceed with the construction of the second Legendrian loop $\Xi_2$, also associated to the Legendrian links $\La(4,4s)$. In conjunction with $\Xi_3$, to be described momentarily, and the Legendrian loop $\Xi_1$ above, $\Xi_1,\Xi_2,\Xi_3$ will be the geometric ingredient for Theorem \ref{thm:main2}.

\begin{remark} The Legendrian loops $\Xi_1,\Xi_2,\Xi_3$ are geometrically constructed to algebraically produce an action of the braid group $B_4$ into $\C[\SM(\La(4,4))]$. Intuitively, $\Xi_1,\Xi_2,\Xi_3$ act respectively as the three Artin generators for $B_4$.\hfill$\Box$
\end{remark}

\subsection{The $\Xi_2$-Loop for $\La(4,4)$}\label{ssec:Xi2} Let us now construct the Legendrian loop $\Xi_2$ for the maximal-tb Legendrian links $\La(4,4s)\sse(\R^3,\xi_\st)$, $s\in\N$. We shall describe it using the same notation as in Subsection \ref{ssec:Xi1} above. The Legendrian loop $\Xi_2$ starts with the braid word $(\sigma_1\sigma_2\sigma_3)^{4(s+1)}$ and it is described by the following sequence of Legendrian isotopies:
\begin{align*}
&(\sigma_1\sigma_2\sigma_3\sigma_1\sigma_2\sigma_3\sigma_1\sigma_2\sigma_3\sigma_1\sigma_2\sigma_3)^{(s+1)} \stackrel{\gamma}{\approx} (\sigma_1\sigma_2{\color{blue}{\underline{\sigma_1\sigma_3}}}\sigma_2\sigma_3\sigma_1\sigma_2\sigma_3\sigma_1\sigma_2\sigma_3)^{(s+1)} \\
&\stackrel{R3^a}{\approx} ({\color{blue}{\underline{\sigma_2\sigma_1\sigma_2}}}\sigma_3 \sigma_2\sigma_3\sigma_1\sigma_2\sigma_3\sigma_1\sigma_2\sigma_3)^{(s+1)} \stackrel{\delta}{\approx} (\sigma_1\sigma_2\sigma_3 \sigma_2\sigma_3\sigma_1\sigma_2\sigma_3\sigma_1\sigma_2\sigma_3{\color{blue}{\underline{\sigma_2}}})^{(s+1)} \\
&\stackrel{\gamma}{\approx} (\sigma_1\sigma_2\sigma_3 \sigma_2{\color{blue}{\underline{\sigma_1\sigma_3}}}\sigma_2\sigma_3\sigma_1\sigma_2\sigma_3\sigma_2)^{(s+1)} \stackrel{R3^d}{\approx} (\sigma_1\sigma_2\sigma_3\sigma_2\sigma_1{\color{blue}{\underline{\sigma_2\sigma_3\sigma_2}}}\sigma_1\sigma_2\sigma_3\sigma_2)^{(s+1)} \quad(\Psi^{(1)}_t \text{until here}) \\
&\stackrel{R3^d}{\approx} (\sigma_1\sigma_2\sigma_3{\color{blue}{\underline{\sigma_1\sigma_2\sigma_1}}}\sigma_3\sigma_2\sigma_1\sigma_2\sigma_3\sigma_2)^{(s+1)} \stackrel{R3^a}{\approx} (\sigma_1\sigma_2\sigma_3\sigma_1\sigma_2\sigma_1\sigma_3\sigma_2\sigma_1{\color{blue}{\underline{\sigma_3\sigma_2\sigma_3}}})^{(s+1)} \\
&\stackrel{\gamma^2}{\approx} (\sigma_1\sigma_2\sigma_3\sigma_1\sigma_2{\color{blue}{\underline{\sigma_3\sigma_1}}}\sigma_2{\color{blue}{\underline{\sigma_3\sigma_1}}}\sigma_2\sigma_3)^{(s+1)}	\quad(\Psi^{(2)}_t \text{until here}).
\end{align*}

In each row, the underlined letters -- emphasized in color blue -- represent those crossings which have been affected when performing the indicated Legendrian isotopy, consisting either of a Reidemeister R3 move, a cyclic shift $\delta$ or a commutation $\gamma$. In the above description of $\Xi_2$, we denote by $\Psi^{(1)}_t$ the Legendrian isotopy consisting of the moves performed in the first six equivalences, and we denote by $\Psi^{(2)}_t$ the Legendrian isotopy consisting of the moves performed in the last four equivalences. The decomposition into the two pieces $\Psi^{(1)}_t$ and $\Psi^{(2)}_t$ will be used in Section \ref{sec:mainproof2}. Note that these $\Psi^{(1)}_t$ and $\Psi^{(2)}_t$ pieces for the Legendrian loop $\Xi_2$ are different from the $\Psi^{(1)}_t$ and $\Psi^{(2)}_t$ pieces for the Legendrian loop $\Xi_1$ in Subsection \ref{ssec:Xi1} above; this repeated notation for the pieces is acceptable because we will only be using these pieces to study $\Xi_1$ or $\Xi_2$ one loop at a time, and thus the notation will be clear by context. Finally, note that the sequence starts and ends with the braid word $(\sigma_1\sigma_2\sigma_3)^{4(s+1)}$, and thus defines a Legendrian loop for $\La(4,4s)$ according to the fronts represented by each braid word.

\begin{definition}\label{def:Xi2Loop} Consider $\La(4,4s)\sse(\R^3,\xi_\st)$, the Legendrian isotopy $\Xi_2$ is the Legendrian loop of $\La(4,4s)$ given by the sequence of Legendrian isotopies above once the zero section $S^1\sse(J^1S^1,\xi_\st)$ is satellited to the standard unknot.\hfill$\Box$
\end{definition}

\subsection{The $\Xi_3$-Loop}\label{ssec:Xi3} We now construct the third Legendrian loop $\Xi_3$ for $\La(4,4s)\sse(\R^3,\xi_\st)$, $s\in\N$. The Legendrian loop $\Xi_3$ starts with the braid word $(\sigma_1\sigma_2\sigma_3)^{4(s+1)}$ and it is described by the following sequence of Legendrian isotopies:
\begin{align*}
&(\sigma_1\sigma_2\sigma_3\sigma_1\sigma_2\sigma_3\sigma_1\sigma_2\sigma_3\sigma_1\sigma_2\sigma_3)^{(s+1)} \stackrel{\gamma}{\approx} (\sigma_1\sigma_2{\color{blue}{\underline{\sigma_1\sigma_3}}}\sigma_2\sigma_3\sigma_1\sigma_2\sigma_3\sigma_1\sigma_2\sigma_3)^{(s+1)} \\
&\stackrel{R3^a}{\approx} ({\color{blue}{\underline{\sigma_2\sigma_1\sigma_2}}}\sigma_3 \sigma_2\sigma_3\sigma_1\sigma_2\sigma_3\sigma_1\sigma_2\sigma_3)^{(s+1)} \stackrel{R3^a}{\approx} (\sigma_2\sigma_1{\color{blue}{\underline{\sigma_3\sigma_2\sigma_3}}}\sigma_3\sigma_1\sigma_2\sigma_3\sigma_1\sigma_2\sigma_3)^{(s+1)}  \quad(\Psi^{(1)}_t \text{until here})\\
&\stackrel{\gamma}{\approx} (\sigma_2{\color{blue}{\underline{\sigma_3\sigma_1}}}\sigma_2\sigma_3\sigma_3\sigma_1\sigma_2\sigma_3\sigma_1\sigma_2\sigma_3)^{(s+1)} \stackrel{\gamma^2}{\approx} (\sigma_2\sigma_3\sigma_1\sigma_2{\color{blue}{\underline{\sigma_1\sigma_3\sigma_3}}}\sigma_2\sigma_3\sigma_1\sigma_2\sigma_3)^{(s+1)} \\
&\stackrel{R3^a}{\approx} (\sigma_2\sigma_3{\color{blue}{\underline{\sigma_2\sigma_1\sigma_2}}}\sigma_3\sigma_3\sigma_2\sigma_3\sigma_1\sigma_2\sigma_3)^{(s+1)} \stackrel{R3^a}{\approx} ({\color{blue}{\underline{\sigma_3\sigma_2\sigma_3}}}\sigma_1\sigma_2\sigma_3\sigma_3\sigma_2\sigma_3\sigma_1\sigma_2\sigma_3)^{(s+1)} \quad(\Psi^{(2)}_t \text{until here}).\\
&\stackrel{\gamma}{\approx} (\sigma_3\sigma_2{\color{blue}{\underline{\sigma_1\sigma_3}}}\sigma_2\sigma_3\sigma_3\sigma_2\sigma_3\sigma_1\sigma_2\sigma_3)^{(s+1)} \stackrel{R3^d}{\approx} (\sigma_3\sigma_2\sigma_1{\color{blue}{\underline{\sigma_2\sigma_3\sigma_2}}}\sigma_3\sigma_2\sigma_3\sigma_1\sigma_2\sigma_3)^{(s+1)} \\
&\stackrel{R3^d}{\approx} (\sigma_3{\color{blue}{\underline{\sigma_1\sigma_2\sigma_1}}}\sigma_3\sigma_2\sigma_3\sigma_2\sigma_3\sigma_1\sigma_2\sigma_3)^{(s+1)} \stackrel{\gamma}{\approx} (\sigma_3\sigma_1\sigma_2{\color{blue}{\underline{\sigma_3\sigma_1}}}\sigma_2\sigma_3\sigma_2\sigma_3\sigma_1\sigma_2\sigma_3)^{(s+1)} \\
&\stackrel{\gamma}{\approx} (\sigma_3\sigma_1\sigma_2\sigma_3\sigma_1\sigma_2\sigma_3\sigma_2{\color{blue}{\underline{\sigma_1\sigma_3}}}\sigma_2\sigma_3)^{(s+1)} \stackrel{R3^d}{\approx} (\sigma_3\sigma_1\sigma_2\sigma_3\sigma_1\sigma_2\sigma_3\sigma_2\sigma_1{\color{blue}{\underline{\sigma_2\sigma_3\sigma_2}}})^{(s+1)} \\
&\stackrel{R3^d}{\approx} (\sigma_3\sigma_1\sigma_2\sigma_3\sigma_1\sigma_2\sigma_3{\color{blue}{\underline{\sigma_1\sigma_2\sigma_1}}}\sigma_3\sigma_2)^{(s+1)} \stackrel{\gamma}{\approx} (\sigma_3\sigma_1\sigma_2\sigma_3\sigma_1\sigma_2\sigma_3\sigma_1\sigma_2{\color{blue}{\underline{\sigma_3\sigma_1}}}\sigma_2)^{(s+1)} \\
&\stackrel{\delta}{\approx} (\sigma_1\sigma_2\sigma_3\sigma_1\sigma_2\sigma_3\sigma_1\sigma_2\sigma_3\sigma_1\sigma_2{\color{blue}{\underline{\sigma_3}}})^{(s+1)}=(\sigma_1\sigma_2\sigma_3)^{4(s+1)} \quad(\Psi^{(3)}_t \text{until here}).
\end{align*}

Note again that the two pieces $\Psi^{(1)}_t,\Psi^{(2)}_t$ for this Legendrian loop $\Xi_3$ differ from the $\Psi^{(1)}_t$ and $\Psi^{(2)}_t$ pieces for the Legendrian loops $\Xi_1,\Xi_2$ in Subsections \ref{ssec:Xi1} and \ref{ssec:Xi2} above.

\begin{definition}\label{def:Xi3Loop} Consider $\La(4,4s)\sse(\R^3,\xi_\st)$, the Legendrian isotopy $\Xi_3$ is the Legendrian loop of $\La(4,4s)$ given by concatenating the sequence of Legendrian isotopies above once the zero section $S^1\sse(J^1S^1,\xi_\st)$ is satellited to the standard unknot.\hfill$\Box$
\end{definition}

The loops $\Sigma_1,\delta^2,\Xi_1,\Xi_2$ and $\Xi_3$ are the needed geometric ingredients in our proof of Theorems \ref{thm:main} and Theorems \ref{thm:main2}. The Legendrian loops $\Sigma_1,\delta^2$ will give rise to the modular action, and $\Xi_1,\Xi_2,\Xi_3$ to the faithful representation of $M_{0,4}$. From a contact topology viewpoint, it is quite outstanding that the infinitely many Lagrangian fillings in Corollary \ref{cor:infinitefillings} can arise in this direct and explicit manner. Let us now move to the algebraic invariants that we shall use in order to build the representations of the modular group $\PSL(2,\Z)$ and the mapping class group $M_{0,4}$.

\begin{remark} The reader is invited to discover the analogue of $\Xi_i$,$1\leq i\leq 3$, for the positive braid $\beta=(\sigma_1\sigma_2\cdot\ldots\cdot\sigma_{n-1})^{n(s+1)}$. These are Legendrian loops for the $n$-component Legendrian links $\La(n,ns)$. We shall nevertheless not need these loops in the present article and thus we do not presently discuss them.\hfill$\Box$
\end{remark}
%
%

\section{Microlocal Legendrian Invariants}\label{sec:algebra}

In this section we introduce the algebraic invariants that we use in order to construct the representations in Theorems \ref{thm:main} and \ref{thm:main2}. These are Legendrian invariants arising from microlocal analysis and the study of constructible sheaves on stratified spaces, as introduced by M. Kashiwara and P. Schapira in the works \cite{GKS_Quantization,KashiwaraSchapira_Book}. The articles \cite{STWZ,STZ} have recently been developing these Legendrian invariants. The present manuscript highlights a remarkable application of these invariants to the study of Lagrangian fillings.\\

Let $\La\sse (\R^3,\xi_\st)$ be a Legendrian link and identify the standard contact 3-space with the positive hemisphere bundle $(T^{\infty,+}(\R^2),\xi_\st)$ of the real 2-plane. Let $\Sh_\La(\R^2,\C)$ be the derived dg-category of constructible sheaves of $\C$-vector spaces on $\R^2$ with singular support intersecting $T^\infty\R^2$ within the Legendrian $\La$. Suppose that $\mbox{rot}(\La)=0$ and consider the microlocal monodromy functor $\mu\mbox{mon}:\Sh_\La(\R^2,\C)\longrightarrow\mbox{Loc}(\La)$ to the category of local systems of complexes of $\C$-vector spaces \cite[Section 5.1]{STZ}. This allows us to consider the following moduli of objects
$$\SM^\circ(\La):=\{\SF^\bullet\in \Sh_\La(\R^2,\C):\mbox{rk}(\mu\mbox{mon}(\SF^{\bullet}))=1,\mu\mbox{mon}(\SF^{\bullet})\mbox{ concentrated in degree }0\}.$$
It is shown in \cite{GKS_Quantization,STZ} that the category $\Sh_\La(\R^2,\C)$, and in particular $\SM^\circ(\La)$, is a Legendrian invariant of $\La$. In the present article, we restrict to Legendrian links $\La\sse(\R^3,\xi_\st)$ which arise as $\La(\beta)$ for a positive braid $\beta$. For this class of Legendrian links, $\mbox{rot}(\La(\beta))=0$ and there exists a binary Maslov potential. Indeed, the braid piece carries the zero Maslov potential and satelliting to the standard Legendrian unknot - with its standard front -  increases the Maslov potential by exactly one.

\subsection{The Brou\'e-Deligne-Michel Description}\label{ssec:openbottsamelson} In order to directly compute with the moduli spaces $\SM^\circ(\La(\beta))$ and construct the representations in Theorems \ref{thm:main} and \ref{thm:main2}, we require a more explicit description of the moduli space $\SM^\circ(\La(\beta))$. This description is available due to the work \cite{STZ}, which proves that $\SM^\circ(\La(\beta))$ is isomorphic to a classical moduli BS($\beta$), modulo the gauge action, associated to a braid $\beta$ by M. Brou\'e-J. Michel \cite{BroueMichel} and Deligne \cite{Deligne}.

Let $G=\GL_k(\C)$ and $B\sse G$ the Borel subgroup of upper triangular matrices. The quotient $G/B$ is the flag variety, whose points parametrize complete flags $V^\bullet$ of vector subspaces of $\C^k$. The Bruhat decomposition
$$\displaystyle G/B=\bigsqcup_{w\in W}BwB/B$$
implies that the relative position of a pair of flags $(V_1^\bullet,V_2^\bullet)$ is determined by an element $s\in S_k=\mbox{Weyl}{(G)}$ of the Weyl group, in this case a permutation in the symmetric group. Consider the Artin generators $\sigma_i\in B_k$, $1\leq i\leq k-1$, and denote by $\overline{\sigma_i}$ the image of $\sigma_i$ under the projection $B_k\lr S_k$ from the braid group to the $A_{k-1}$-Coxeter group $S_k$. Given a flag $V^\bullet$ and a permutation $s\in S_k$, let $S_s(V^\bullet)$ be the set of flags in relative $s$-position with respect to $V^\bullet$.

\begin{definition}
Let $\beta$ be a positive braid word
$$\beta=\prod_{j=1}^{l(\beta)}\sigma_{i_j},\quad 1\leq i_j\leq k-1,$$
and consider the subset
$$\mbox{BS}(\beta):=\{(V^\bullet_1,\ldots,V^\bullet_{l(\beta)})\in (G/B)^{l(\beta)}:V^\bullet_{m+1}\in S_{\overline{\sigma}_{i_{m}}}(V^\bullet_{m}),1\leq m\leq l(\beta)\}\sse(G/B)^{l(\beta)},$$
where the index $1\leq m\leq l(\beta)$ is understood cyclically modulo $l(\beta)$, i.e. the condition for $m=l(\beta)$ reads $V^\bullet_1\in S_{\overline{\sigma}_{i_{l(\beta)}}}(V^\bullet_{l(\beta)})$. By definition, BS($\beta$) is said to be the open Bott-Samelson variety associated to $\beta$.\hfill$\Box$
\end{definition}
%


For each $\beta$, the group $G$ acts on the open Bott-Samelson variety BS($\beta$) diagonally on the left, given that the flag variety $G/B$ is given by the $B$-action on the right. The article \cite[Section 6]{STZ} identifies $\SM^\circ(\Lambda(\beta))$ with the quotient $G\backslash$BS($\beta$). It is a consequence of this identification that our moduli space $\SM^\circ(\Lambda(\beta))$ can be described as follows.

Choose a set of points $\{\theta_0,\theta_1,\ldots,\theta_{l(\beta)}\}\in S^1$ such that the vertical lines $\{\theta_m\}\times\R$, $0\leq m\leq l(\beta)$ do not intersect the front $\beta\sse S^1\times\R$ at a crossing and there exists a unique crossing of $\beta$ between $\{\theta_m\}\times\R$ and $\{\theta_{m+1}\}\times\R$, $0\leq i\leq l(\beta)$. Then $\SM^\circ(\Lambda(\beta))$ is the moduli space given by associating a complete flag $V^\bullet_m$ along each vertical line $\{\theta_m\}\times\R$ such that $V^\bullet_0=V^\bullet_{l(\beta)}$ and two flags $V^\bullet_m$ and $V^\bullet_{m+1}$ differ only and exactly in their $i_{m+1}$-dimensional subspaces for all $1\leq m\leq l(\beta)$, modulo the gauge group action of $\GL_k(\C)$. This description in terms of BS($\beta$) will be used in Sections \ref{sec:mainproof1} and \ref{sec:mainproof2}.

\subsection{Moduli of Framed Sheaves} In the proof of Theorems \ref{thm:main} and \ref{thm:main2} we shall need a framed enhancement $\SM(\La(\beta),\tau)$ of the Bott-Samelson varieties $\SM^\circ(\La(\beta))$. In precise terms, the points of $\SM(\La(\beta),\tau)$ are given by the $l(\beta)$-tuples of flags $[(V^\bullet_1,\ldots,V^\bullet_{l(\beta)})]\in\SM^\circ(\La(\beta))$ equipped with trivializations $\tau$ for the stalks at a specified set of points. In this case, we choose the set of points such that the set contains exactly one point for each region where the constructible sheaf has a 1-dimensional stalk. Given that they are in bijection, we will interchangeably speak of these points or the open strata in the front diagram that contain them, these open strata shall also be referred to as regions. Hence, in the language of Bott-Samelson varieties, the trivialization $\tau$ consists of a series of isomorphisms
$$V^{(1)}_m\cong\C,\quad1\leq m\leq l(\beta).$$
In our context, the moduli spaces of framed sheaves $\SM(\La(\beta),\tau)$ are algebraic varieties \cite{STZ}. It should be emphasized that the moduli space $\SM(\La(\beta),\tau)$ depends on the choice of trivialization $\tau$. In our choice above,  $\SM(\La(\beta),\tau)$ shall depend on the choice of braid word $\beta$. Indeed, the length of the tuple is precisely $l(\beta)$. Nevertheless, the article \cite{STWZ} shows that a Legendrian isotopy generates an equivalence of moduli space $\SM(\La(\beta),\tau)$ of framed sheaves, with the trivialization, and its region, being pushed forward under the isotopy. Thus, in studying the action of a Legendrian loop on $\SM(\La(\beta),\tau)$ we identify the moduli spaces of framed sheaves along the Legendrian isotopy and compare the action at the canonically identified endpoints of the Legendrian loop.

Explicitly, let $\{\Psi_t\}_{t\in[0,1]}$ be a Legendrian loop based at the identity, i.e. $\Psi_0=\mbox{Id}$. By \cite[Section 2]{STWZ}, there is a canonical isomorphism between the moduli spaces $\SM(\La,\tau)$ and $\SM(\Psi_t(\La),(\Psi_t)_*\tau)$ for all $t\in[0,1]$. By virtue of being a Legendrian loop, $\Psi_1=\mbox{Id}$ and thus we obtain an algebraic automorphism $\Psi\in\Aut(\SM(\La,\tau))$ of the moduli space $\SM(\La,\tau)$. This automorphism is to be understood as the monodromy of the Legendrian loop $\{\Psi_t\}_{t\in[0,1]}$, in line with T. K\'alm\'an's \cite[Section 3]{Kalman} monodromy invariant. The automorphism $\Psi\in\Aut(\SM(\La,\tau))$ in turn induces an automorphism $\Psi^*\in\Aut(\C[\SM(\La,\tau)])$ in the coordinate ring of regular functions on $\SM(\La,\tau)$.

In Theorems \ref{thm:main} and \ref{thm:main2} the focus will be on two moduli spaces $\SM(\La(\beta),\tau)$ for the two braid words $\beta=(\sigma_1\sigma_2)^9$ and $\beta=(\sigma_1\sigma_2\sigma_3)^{8}$ and a chosen trivialization $\tau$.

\subsection{Ingredients on $\SL_3$-webs}\label{ssec:SL3webs} The argument for the faithfulness in the statement of Theorem \ref{thm:main}, as presented in Section \ref{sec:mainproof1}, requires the study of the coordinate ring $\C[\Gr(3,9)]$. We need regular functions beyond the Pl\"ucker coordinates in $\C[\Gr(3,9)]$ because the pull-back of some of the Pl\"ucker coordinates under the (action on certain moduli spaces induced by our) Legendrian loops are no longer Pl\"ucker coordinates. Thus, we provide in this subsection the ingredients that we use to study $\C[\Gr(3,9)]$. They were developed in \cite{KuperbergSpiders} originally, and we will use the notation and perspective established in \cite{FominPylyavskyy}.

Consider a closed disk $D$ with $m$ marked points on the boundary. By definition, a {\it tensor diagram} in $D$ for $\SL_3$ is a finite bipartite graph drawn in $D$ with a bipartition of its vertex set into black and white color sets such that:
\begin{itemize}
	\item[-] The boundary marked points of $D$ are black vertices of the graph, and they are the only vertices of the graph at the boundary.
	\item[-] The vertices which are not marked points, in the interior, are trivalent.
\end{itemize}
The case of interest in this manuscript is $m=9$ marked points at the boundary.

\begin{center}
	\begin{figure}[h!]
		\centering
		\includegraphics[scale=0.7]{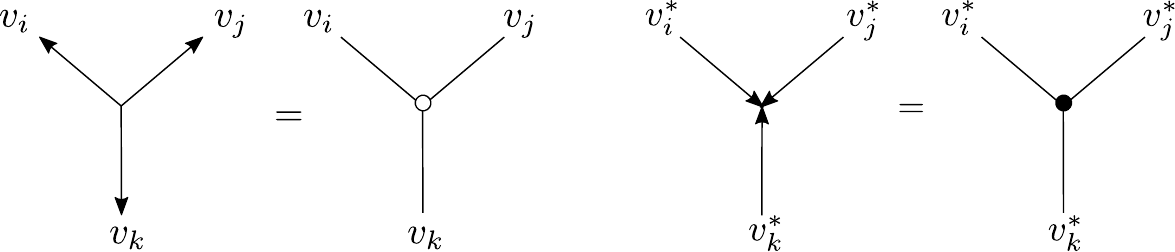}
		\caption{The diagrams on each of the left hand sides of the equalities are called $A_2$-spiders in \cite[Section 4]{KuperbergSpiders}. The diagrams on the right hand sides follow the notational convention of \cite{FominPylyavskyy}, white vertices are sources and black vertices are sinks. The white tripod represents $\det(v_iv_jv_k)$, $v_i,v_j,v_k\in V\cong\C^3$, which is the $\SL_3$-invariant tensor given by a fixed volume form $V^{\otimes 3}\longrightarrow \C$, and the black tripod represents its dual.}
		\label{fig:TensorDiagramsSet4}
	\end{figure}
\end{center}

Let $V=\C^3$ be a vector space endowed with a volume form. Suppose we assign a vector $v\in V$ to each black vertex, and a covector $v^*\in V^*$ to each white vertex. Two basic $\SL_3$-invariant tensors associated to $V$ are the volume form $V^{\otimes 3}\longrightarrow \C$ and the dual form $(V^*)^{\otimes 3}\longrightarrow \C$. For the purposes of this manuscript, they are diagrammatically encoded by a white tripod and a black tripod, respectively, as depicted in Figure \ref{fig:TensorDiagramsSet4}. This follows the notation of \cite{FominPylyavskyy}, with white and black vertices, but note that these diagrammatics were previously studied in \cite{KuperbergSpiders} for rank 2 algebras; in particular, $\SL_3$ is associated to the $A_2$-Dynkin diagram, and these tensor diagrams were called $A_2$-spiders by G. Kuperberg. The canonical pairing $V\otimes V^*\longrightarrow\C$ is diagrammatically given by an edge between a black and white vertex, i.e. an edge can also be considered as the identity in $V$ if we identified $V\cong V^*$.

\begin{center}
	\begin{figure}[h!]
		\centering
		\includegraphics[scale=0.9]{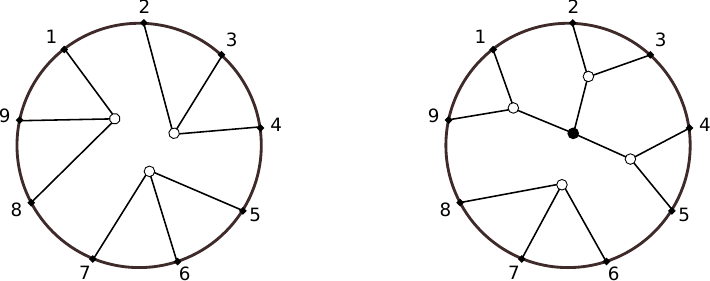}
		\caption{Two examples of $\SL_3$-webs for the Grassmannian $\Gr(3,9)$. As invariant functions, and as elements of the coordinate ring $\C[\Gr(3,9)]$, these tensor diagrams represent the functions $\det(v_2v_3v_4)\det(v_5v_6v_7)\det(v_8v_9v_1)$, for the diagram on the left, and $\det(v_9\times v_1,v_2\times v_3,v_4\times v_5)\det(v_6v_7v_8)$ for the diagram on the right.}
		\label{fig:TensorDiagramsSet3}
	\end{figure}
\end{center}

Now, suppose that vectors $v_1,\ldots,v_m\in V$ are assigned to the $m$ marked points at the boundary of $D$, one vector per marked point. Then, a tensor diagram can be used to define a $\C$-scalar by repeated contraction using the basic $\SL_3$-invariant tensors. For instance, Figure \ref{fig:TensorDiagramsSet3} give two examples of tensor diagrams and their associated functions for $m=9$, and see \cite{FominPylyavskyy} and \cite[Section 9]{Fraser} for more details.

Finally, a point in the (affine cone of the) Grassmannian $\Gr(3,9)$ will be represented by an ordered tuple of $9$ vectors in $V$, modulo the appropriate action. In this manner, a tensor diagram gives rise to a regular function in the coordinate ring $\C[\Gr(3,9)]$. For instance, Figure \ref{fig:TensorDiagramsSet3} (left) represents the product $P_{234}P_{567}P_{891}\in\C[\Gr(3,9)]$, where $P_{ijk}=v_i\wedge v_j\wedge v_k$ is a Pl\"ucker coordinate. See Subsection \ref{ssec:FaithfulPSL} for further examples.

\begin{remark}
	We conclude with a piece of terminology. A planar tensor diagram is often called a {\it web} in the literature. This is the reason that this diagrammatic calculus is referred to as web combinatorics, and we refer to the webs associated planar tensor diagrams for $\SL_3$ as {\it $\SL_3$-webs}, as in \cite{Fraser}. Following \cite{KuperbergSpiders}, a web is non-elliptic if it contains no 2-cycles based at a boundary vertex, and if all of its faces formed
	by interior vertices are bounded by at least six sides. G. Kuperberg showed in \cite{KuperbergSpiders} that (non-elliptic) webs can be used to construct bases for many rings of $\SL_3$-invariants.\hfill$\Box$
\end{remark} 

In the next two sections we prove Theorem \ref{thm:main} and \ref{thm:main2}. These two proofs are independent of each other. The reader is nevertheless encouraged to read the proof of Theorem \ref{thm:main} first, as it also sets the main techniques and notations for the proof of Theorem \ref{thm:main2}.

\section{The $\PSL(2,\Z)$ representation for $\La(3,6)$}\label{sec:mainproof1}

\begin{center}
	\begin{figure}[h!]
		\centering
		\includegraphics[scale=0.7]{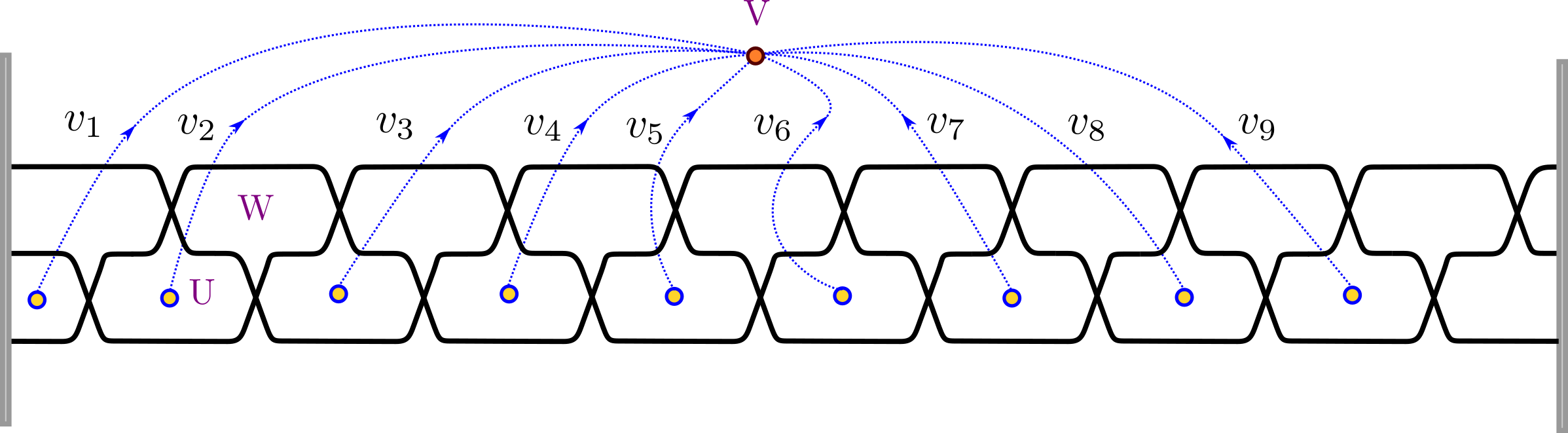}
		\caption{The identification of $\SM(\La(3,6))$ with the positroid cell in $\Gr(3,9)$.}
		\label{fig:T39SheavesToGrassmannian}
	\end{figure}
\end{center}

Let us prove Theorem \ref{thm:main}. For that, we shall compute the action of the two Legendrian loops $\Sigma_1,\delta^2$ constructed in Section \ref{sec:loops} into the coordinate ring $\C[\SM(\La(\beta),\tau)]$ of the framed Bott-Samelson variety $\SM(\La(\beta),\tau)$, where the braid is fixed to be $\beta=(\sigma_1\sigma_2)^9$ and the trivialization $\tau$ is given at the 1-dimensional stalks depicted as dots Figure \ref{fig:T39SheavesToGrassmannian}, where the braid $\beta$ for the Legendrian link $\La(3,6)$ is also depicted. These monodromy invariants $\Sigma_1^*,(\delta^2)^*\in\Aut(\C[\SM(\La(\beta),\tau)])$ will be shown to be non-trivial and generate an action of an infinite group on the coordinate ring $\C[\SM(\La(\beta),\tau)]$.
\subsection{The Monodromy Effect within $\Gr(3,9)$}\label{ssec:GrassmannianEmbedding1} The first step of the argument is to identify the moduli $\SM(\La(\beta),\tau)$ with the positroid stratum $\Pi_r\sse\Gr(3,9)$ in the projective Grassmannian $\Gr(3,9)$, where $r$ is the cyclic rank matrix associated to the positive braid $(\sigma_1\sigma_2)^9$. The canonical embedding of $\SM(\La(\beta),\tau)$ into $\Gr(3,9)$, with image $\Pi_r$, is obtained as follows \cite[Section 3.2]{STWZ}. Given a point $(V^\bullet_1,\ldots,V^\bullet_{18})\in\SM(\La(\beta),\tau)$, the $9$-tuple of vectors $$(v_1,v_2,\ldots,v_9)\in(V^{(1)}_1,V^{(1)}_3,V^{(1)}_5,\ldots,V^{(1)}_{17}),$$ modulo the $\GL_9(\C)$-action, defines a point in $\Gr(3,9)$, where the choice of vectors is given by the framing. These vectors $(v_1,v_2,\ldots,v_9)$ are depicted in Figure \ref{fig:T39SheavesToGrassmannian}. The advantage of this algebraic embedding $\SM(\La(\beta),\tau)\lr\Gr(3,9)$ is that it allows us to use elements in the homogeneous coordinate ring of $\Gr(3,9)$ restricted to $\SM(\La(\beta),\tau)$ in order to study the effect of the monodromies $\Sigma_1^*,(\delta^2)^*\in\Aut(\C[\SM(\La(\beta),\tau)])$. We shall henceforth denote the framed moduli space by $\SM(\La(\beta))$, where the trivialization $\tau$ is implicitly chosen to be as above.

\begin{remark} Consider three vector spaces $U,W,V$ of dimensions $\dim(U)=1$, $\dim(W)=2$ and $\dim(V)=3$. A framed constructible sheaf $\SF\in\SM(\La(\beta),\tau)$ has stalks isomorphic to $U$, $W$ and $V$ as depicted in Figure \ref{fig:T39SheavesToGrassmannian}. The $9$-tuple of vectors described above can also be obtained by parallel transport of the stalk of $\SF\in\SM(\La(\beta))$ in the $U$-region to the $V$-region along the dashed paths depicted in Figure \ref{fig:T39SheavesToGrassmannian}. Note that it does not matter whether a dashed arrow passes a crossing from its left or its right.\hfill$\Box$
\end{remark}

Let us now analyze the action of the Legendrian loops $\Sigma_1$ and $\delta^2$ on the coordinate ring of $\SM(\La(\beta))$ by studying their action on the $9$-tuples of vectors $(v_1,v_2,\ldots,v_9)$. For that, we must identify the explicit effect of each of the Legendrian isotopies constituting $\Sigma_1$ and $\delta^2$. These consist of cyclic shifts and Reidemeister III moves.

The effect of the Legendrian isotopy $\delta^2$ described in Subsection \ref{ssec:Legloop}.(ii) above is precisely the cyclic shift on the $9$-tuple of vectors:
$$\delta^2(v_1,v_2,\ldots,v_9)=(v_9,v_2,\ldots,v_1).$$

\begin{center}
	\begin{figure}[h!]
		\centering
		\includegraphics[scale=0.8]{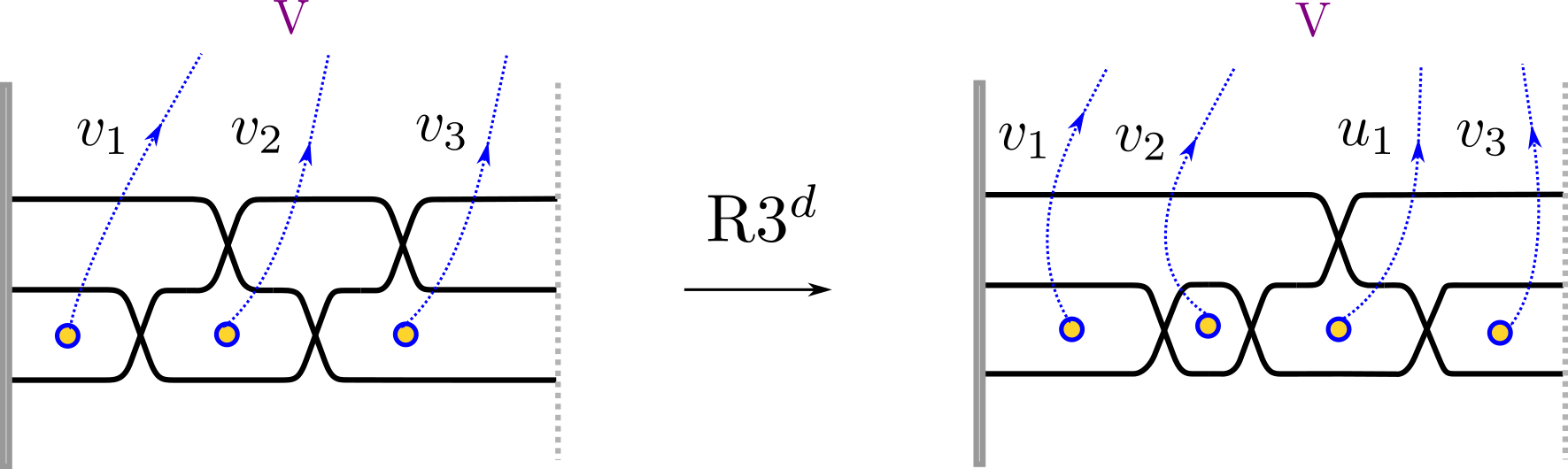}
		\caption{The effect of a Reidemeister R3$^d$ move on $\SM(\La(\beta))$.}
		\label{fig:T39Sheaves_FirstR3}
	\end{figure}
\end{center}

The effect of Reidemeister III moves is more interesting. Indeed, the Reidemeister R3$^{d}$ introduces a $U$-region and thus contributes to a vector $u_1$, whereas the Reidemeister R3$^{a}$, conversely, reduces the number of $U$-regions by exactly one, thus making a vector disappear. Figure \ref{fig:T39Sheaves_FirstR3} depicts the case where the $3$-tuple $(v_1,v_2,v_3)$, in the region given by the braid $\sigma_1\sigma_2\sigma_1\sigma_2$ becomes the $4$-tuple $(v_1,v_2,u_1,v_3)$ for the braid $\sigma_1\sigma_1\sigma_2\sigma_1$.

In terms of the $5$-tuple of flags $(V^\bullet_1,\ldots,V^\bullet_5)$, associated to the braid $\sigma_1\sigma_2\sigma_1\sigma_2$, as described in Section \ref{sec:algebra}, the vectors are $v_1=V^{(1)}_1$, $v_2=V^{(1)}_2=V^{(1)}_3$ and $v_3=V^{(1)}_4=V^{(1)}_5$. Performing the descending Reidemeister III in Figure \ref{fig:T39Sheaves_FirstR3} yields the new vector $u_1$, whose direction is uniquely defined by $V^{(2)}_1\cap V^{(2)}_5$ and the normalization is given by the framing. The following proposition describes the algebraic effect of $\Sigma_1$:


\begin{prop}\label{prop:EffectSigma1} The Legendrian loop $\Sigma_1$ induces the morphism
$$(v_1,v_2,v_3;v_4,v_5,v_6;v_7,v_8,v_9)\longmapsto(v_2,u_1,v_3;v_5,u_2,v_6;v_8,u_3,v_9),$$
where $u_1,u_2$ and $u_3$ are given by the intersections
$$u_1\in\langle v_1,v_2\rangle\cap \langle v_3,v_4\rangle,\quad u_2\in\langle v_4,v_5\rangle\cap \langle v_6,v_7\rangle,\quad u_3\in\langle v_7,v_8\rangle\cap \langle v_9,v_1\rangle,$$
and the three normalizing conditions $v_1v_2=v_2u_1$, $v_4v_5=v_5u_2$ and $v_7v_8=v_8u_3$ in $V\wedge V$.
\hfill$\Box$
\end{prop}

\begin{center}
	\begin{figure}[h!]
		\centering
		\includegraphics[scale=0.8]{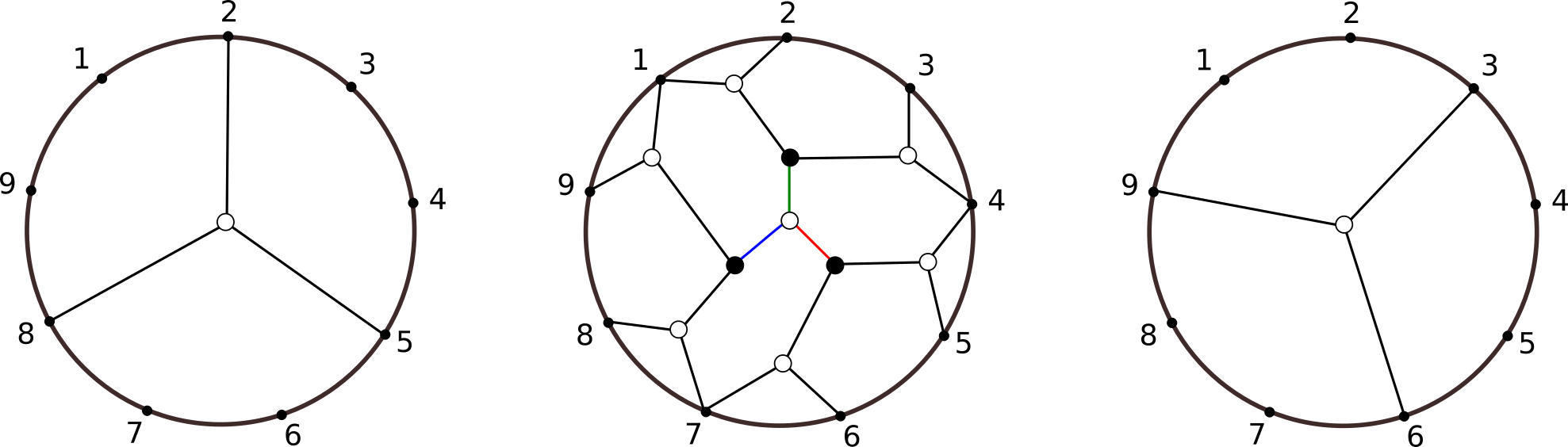}
		\caption{The webs associated to $(\Sigma_1)^*P_{147}$, on the left, $(\Sigma_1)^*P_{258}$, in the center, and $(\Sigma_1)^*P_{369}$, on the right, according to Proposition \ref{prop:EffectSigma1}. In the web for $(\Sigma_1)^*P_{258}$ we have depicted $u_1$ in green, $u_2$ in red and $u_3$ in blue.}
		\label{fig:EffectXi1Plucker}
	\end{figure}
\end{center}

Figure \ref{fig:EffectXi1Plucker} depicts instances of the conclusion of Proposition \ref{prop:EffectSigma1} in terms of G. Kuperberg's $\SL_3$-web combinatorics \cite[Section 4]{KuperbergSpiders}, see Subsection \ref{ssec:SL3webs} above. The reader is also referred to \cite{FominPylyavskyy,KuperbergKhovanov} for the basics of planar tensor diagrams and $\SL_3$-webs, which we shall use in Subsection \ref{ssec:FaithfulPSL}. In particular, Figure \ref{fig:EffectXi1Plucker} displays the pull-backs $(\Sigma_1)^*P_{147}$, $(\Sigma_1)^*P_{258}$ and $(\Sigma_1)^*P_{369}$ of three Pl\"ucker coordinates $P_{ijk}\in\C[\Gr(3,9)]$, $1\leq i<j<k\leq 9$, where $P_{ijk}=v_i\wedge v_j\wedge v_k$. In particular, Proposition \ref{prop:EffectSigma1} implies $(\Sigma_1)^*P_{147}=P_{258}$ and $(\Sigma_1)^*P_{369}=P_{369}$.

We shall provide the proof of Proposition \ref{prop:EffectSigma1} momentarily. Let us however first conclude the proof of Theorem \ref{thm:main} assuming Proposition \ref{prop:EffectSigma1}.

\subsection{The Faithful $\PSL(2,\Z)$-Action}\label{ssec:FaithfulPSL} For that, we study the monodromy action of the subgroup $\Gamma=\langle[\Sigma_1],[\delta^2]\rangle\sse\pi_1(\cL(3,6))$ generated by the homotopy classes of the two Legendrian loops $\Sigma_1,\delta^2$ into the set of $9$-tuples of vectors in $\C^3$. In order to show that this action is indeed non-trivial we choose a function $\Delta\in\C[\Gr(3,9)]$ and ensure that the pull-backs of this function are distinct. For our braid $\beta=(\sigma_1\sigma_2)^9$, let us choose the Pl\"ucker coordinate $\Delta=P_{147}$ in $\C[\Gr(3,9)]$, given by $P_{147}(v_1,\ldots,v_9)=v_1\wedge v_4\wedge v_7$. The algebraic claim that needs to be proven is that the monodromy of $\Sigma_1,\delta^2$ induces a faithful $\PSL(2,\Z)$-action on the orbit $\cO(P_{147})$.

\begin{center}
	\begin{figure}[h!]
		\centering
		\includegraphics[scale=0.7]{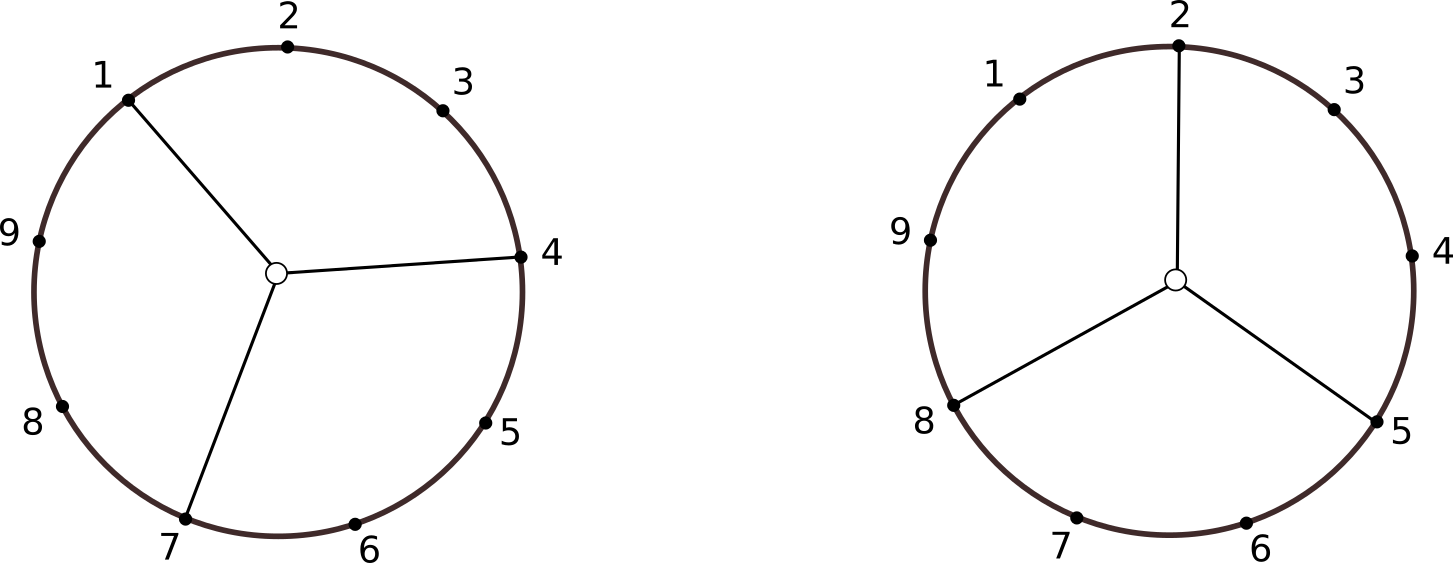}
		\caption{The webs associated to the Pl\"ucker coordinates $P_{147}$, on the left, and $P_{258}$, on the right. In general, $\delta^2$ acts by clockwise rotation on webs.}
		\label{fig:Plucker147}
	\end{figure}
\end{center}

First, let $A=\delta^2$ and $B=\Sigma_1\circ\delta^2$. We have that $A^*(P_{147})=P_{258}$, $(A^2)^*(P_{147})=P_{369}$ and $(A^3)^*(P_{147})=P_{471}=P_{147}$ and thus $A$ generates a $\Z_3$-action on the orbit $\cO(P_{147})$. In general, the action of $A$ and $B$ cannot be exclusively written in terms of Pl\"ucker coordinates. and in order to study our monodromy action we shall be using $\SL_3$-webs, see Subsection \ref{ssec:SL3webs} above and references therein. In terms of $\SL_3$-webs, the diagrams associated to the Pl\"ucker coordinates $P_{147}$ and $A^*(P_{147})=P_{258}$ are depicted in Figure \ref{fig:Plucker147}.

The monodromy of the Legendrian loop $B$ generates a $\Z_2$-action on the orbit $\cO(P_{147})$. Indeed, the square $B^2=\Sigma_1\circ\delta^2\circ \Sigma_1\circ\delta^2$ pull-backs $P_{147}$ as follows
\begin{align*}
(B^*)^2(P_{147})&=(\delta^2)^*\circ(\Sigma_1)^*\circ(\delta^2)^*\circ \Sigma_1^*(P_{147})=(\delta^2)^*\circ(\Sigma_1)^*\circ(\delta^2)^*(P_{258})=\\
&=(\delta^2)^*\circ(\Sigma_1)^*(P_{369})=(\delta^2)^*(P_{369})=P_{147},
\end{align*}
where we have used $(\Sigma_1)^*(P_{369})=P_{369}$, as implied by Proposition \ref{prop:EffectSigma1}. Thus, $A^*$ generates a $\Z_3$-action and $B^*$ generates a $\Z_2$-action. Since the modular group $\PSL(2,\Z)\cong\Z_3*\Z_2$ is a free product, it suffices to show that $A$ and $B$ generate a faithful action with no relations in the subgroup $\langle[A],[B]\rangle$. Following \cite[Section 10]{Fraser}, we will prove this by using the Ping-Pong Lemma \cite[Section II.B]{PingPong1}:

\begin{lemma}[\cite{PingPong1,PingPong2}]\label{lem:pingpong} Let $\Gamma$ be a group acting on a set $X$, let $\Gamma_1,\Gamma_2$ be two subgroups of $\Gamma$, and let $G$ be the subgroup of $\Gamma$ generated by $\Gamma_1$ and $\Gamma_2$. Suppose that $|\Gamma_1|\geq3$ and $|\Gamma_2|\geq2$.
	
Assume that there exist two non-empty subsets $X_1,X_2$ in $X$, with $X_2$ not included in $X_1$, such that
$$\gamma(X_2)\sse X_1,\quad \forall\gamma\in\Gamma_1,\gamma\neq1$$
$$\gamma(X_1)\sse X_2,\quad \forall\gamma\in\Gamma_2,\gamma\neq1.$$
Then $G$ is isomorphic to the free product $\Gamma_1*\Gamma_2$.\hfill$\Box$
\end{lemma}

We apply Lemma \ref{lem:pingpong} for $\Gamma=\pi_1(\cL(3,6))$, $G=\langle[A],[B]\rangle$, $\Gamma_1=\langle[A]\rangle$ and $\Gamma_2=\langle[B]\rangle$, which indeed satisfy $|\Gamma_1|\geq3$ and $|\Gamma_2|\geq2$. The action of $G$ in $X$ is given by the induced monodromy, as described in Section \ref{sec:algebra}. Consider $X=\cO(P_{147})$ to be the orbit of $P_{147}$. Let us now define the Ping-Pong sets $X_1$ and $X_2$. This shall be done in terms of their web diagrams, as follows.

\begin{definition}\label{def:PingPong} The set $X_1\sse\cO(P_{147})$ is the set of all (non-elliptic) webs in $\C[\Gr(3,9)]$ which do not contain any of the pieces in Figure \ref{fig:TensorDiagramsSet1}. That is, a web is in $X\setminus X_1$ if it contains at least one of the pieces in Figure \ref{fig:TensorDiagramsSet1}.
		
		\begin{center}
			\begin{figure}[h!]
				\centering
				\includegraphics[scale=0.7]{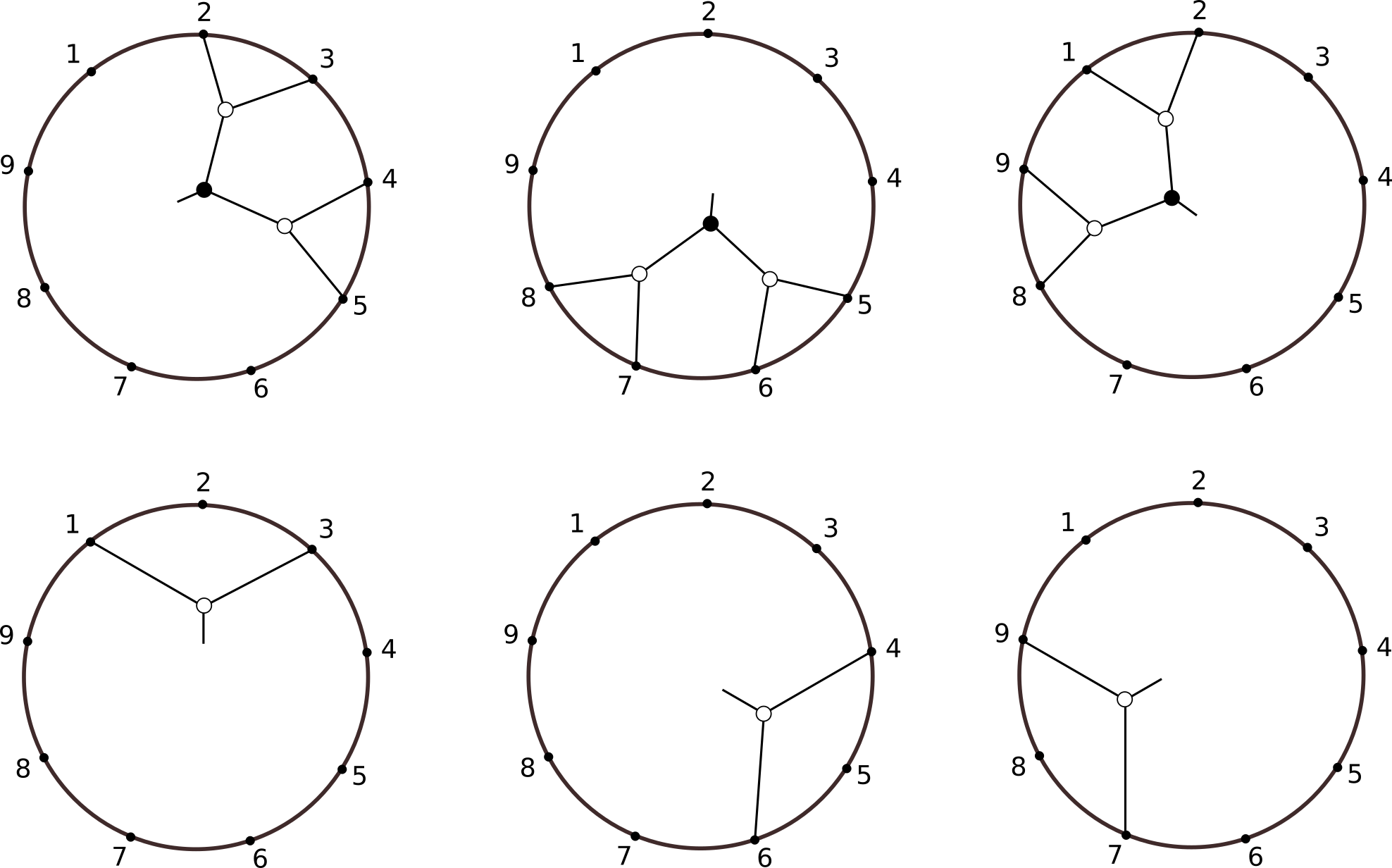}
				\caption{The webs in the set $X_1$ do {\it not} contain any of these six pieces.}
				\label{fig:TensorDiagramsSet1}
			\end{figure}
		\end{center}
	
	Similarly, the set $X_2\sse\cO(P_{147})$ is the set of all (non-elliptic) webs in $\C[\Gr(3,9)]$ which do not contain any of the pieces in Figure \ref{fig:TensorDiagramsSet2}.
	
	\begin{center}
		\begin{figure}[h!]
			\centering
			\includegraphics[scale=0.7]{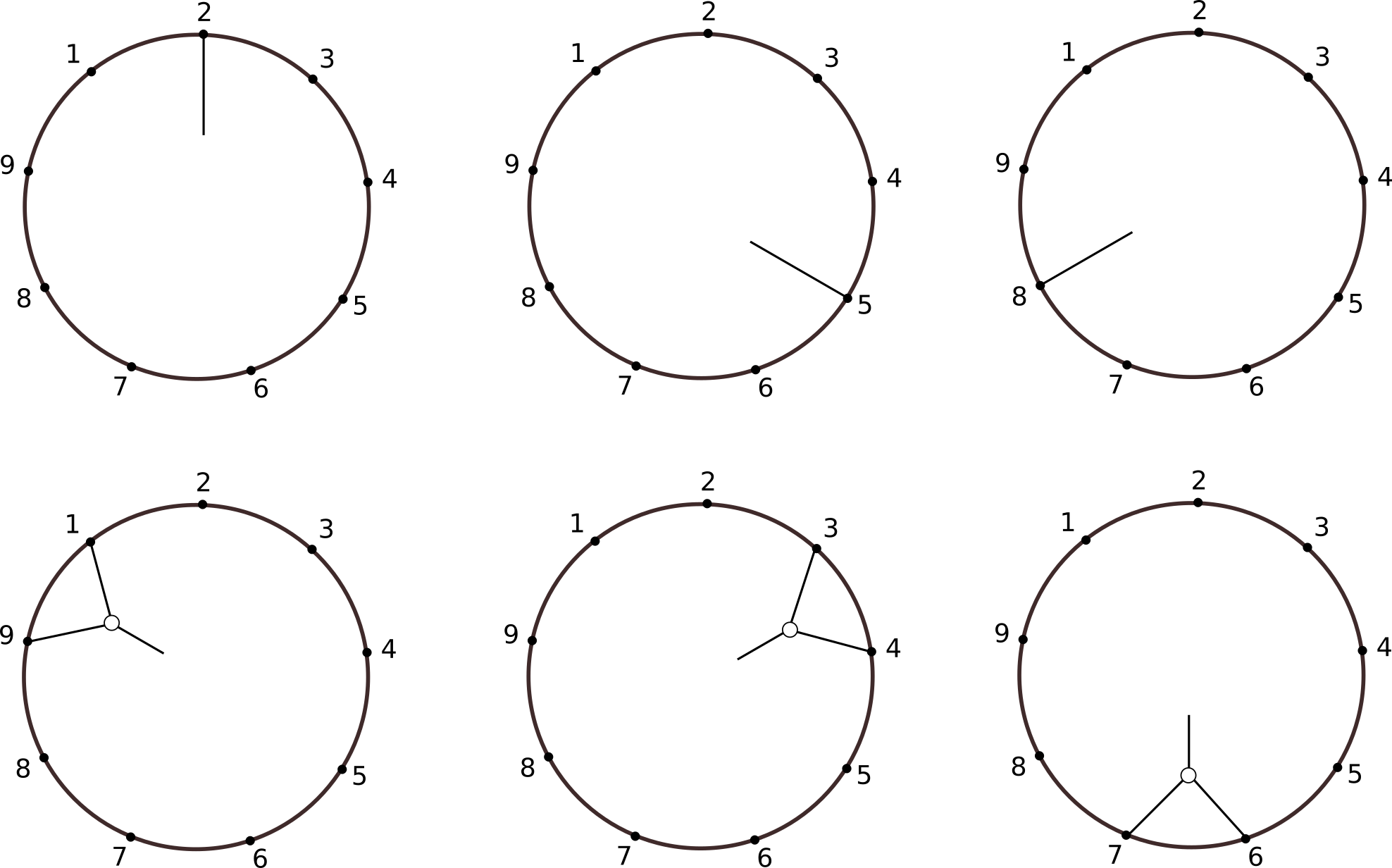}
			\caption{The webs in the set $X_2$ do {\it not} contain any of these six pieces.}
			\label{fig:TensorDiagramsSet2}
		\end{figure}
	\end{center}
\end{definition}

It suffices to prove that $X_1,X_2$ in Definition \ref{def:PingPong} are Ping-Pong sets for the monodromy action. It is useful to remind ourselves that the pull-back $A^*=(\delta^2)^*$ acts by clockwise rotation by $2\pi/9$-radians on the web diagram.

First, let us prove the inclusion $A(X_2)\sse X_1$. Suppose that we have a web $W_2\in X_2$, we need to argue that $A(W_2)$ contains none of the six patterns displayed in Figure \ref{fig:TensorDiagramsSet1}. Suppose $A(W_2)$ contained any of these six, then the clockwise rotation by $4\pi/9$-radians of the diagram $A(W_2)$ will contain a spike at one of the boundary vertices $2$,$5$ or $8$, and thus not be in $X_2$. This rotation by $4\pi/9$-radians of the diagram $A(W_2)$  represents $A^2(A(W_2))=W_2$ since $A^3=\mbox{id}$ in $\cO(P_{147})$, which contradicts $W_2\in X_2$. This shows $A(X_2)\sse X_1$.

Second, let us prove the inclusion $B(X_1)\sse X_2$. Consider a web $W_1\in X_1$, we need to argue that $B(W_1)$ contains none of the six patterns displayed in Figure \ref{fig:TensorDiagramsSet2}. Suppose $B(W_1)$ contained any of the three patterns displayed in the first row of Figure \ref{fig:TensorDiagramsSet2}, i.e. a spike at either one of the boundary vertices $2$,$5$ or $8$. By Proposition \ref{prop:EffectSigma1}, the web $\Sigma_1^*(B(W_1))$ contains one of the three patterns in the first row of Figure \ref{fig:TensorDiagramsSet1} rotated counter-clockwise by an angle of $2\pi/9$-radians. In consequence, the $2\pi/9$-clockwise rotation
$$(\delta^2)^*\circ\Sigma_1^*(B(W_1))=B(B(W_1))=W_1\not\in X_1,$$
of $\Sigma_1^*(B(W_1))$ does not belong to $X_1$, which contradicts $W_1\in X_1$. Thus $B(W_1)$ does not contain any of the three patterns displayed in the first row of Figure \ref{fig:TensorDiagramsSet2}.

Now, suppose that $B(W_1)$ contained any of the three patterns displayed in the second row of Figure \ref{fig:TensorDiagramsSet2}. Proposition \ref{prop:EffectSigma1} implies that the web $\Sigma_1^*(B(W_1))$ contains a counter-clockwise rotated copy, by an angle of $2\pi/9$-radians, of one of the three patterns in the second row of Figure \ref{fig:TensorDiagramsSet1}. Thus, the $2\pi/9$-clockwise rotation $(\delta^2)^*\circ\Sigma_1^*(B(W_1))$ does not belong to $X_2$. This is a contradiction with
$$(\delta^2)^*\circ\Sigma_1^*(B(W_1))=B(B(W_1))=W_1\not\in X_1.$$
Hence $B(W_1)$ cannot contain any of the three patterns displayed in the second row of Figure \ref{fig:TensorDiagramsSet2}. This shows $B(X_1)\sse X_2$, as desired. In conclusion, $X_1$ and $X_2$ are Ping-Pong sets and Lemma \ref{lem:pingpong} implies that $G=\langle[A],[B]\rangle$ is isomorphic to $\PSL(2,\Z)$ and thus the restriction of the monodromy action to this subgroup is a faithful $\PSL(2,\Z)$-representation along the orbit $X=\cO(P_{147})$. This concludes the proof of Theorem \ref{thm:main} once Proposition \ref{prop:EffectSigma1} has been proven.\hfill$\Box$

\subsection{Proof of Proposition \ref{prop:EffectSigma1}} 
\begin{center}
		\begin{figure}[h]
	\centering
	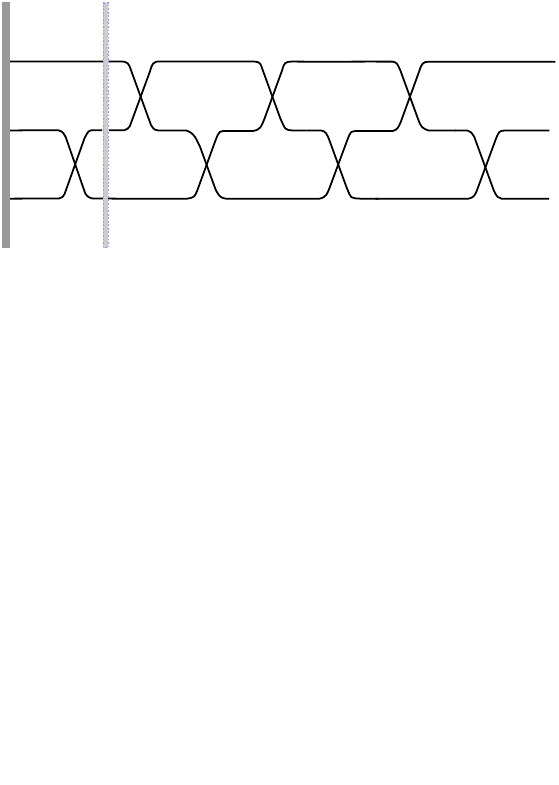
	\caption[]{First window $\beta_0$ of the Legendrian braid $\beta=(\sigma_1\sigma_2)^9$ associated to the Legendrian link $\Lambda(3,6)$. There are total of three windows.}
	\label{fig: T39window}
\end{figure}
\end{center}

Let us consider the braid word $\beta=(\sigma_1\sigma_2)^9$ and consider the braid word given by the piece $\beta_0=(\sigma_1\sigma_2)^3$, such that $\beta=\beta_0^3$ is a concatenation of $\beta_0$ three times. We refer to the piece $\beta_0$ as a window for the braid $\beta$, such that $\beta$ consists of three windows. The Legendrian loop $\Sigma_1$ consists of a cyclic permutation and a sequence of braid equivalences given by the Reidemeister III moves. The braid equivalence can be performed equivariantly over each of the three windows, and hence the morphism induced from $\Sigma_1$ is periodic with respect to this prescribed window decomposition once the shift is applied. It thus suffices to work with one window to describe the morphism. Figure \ref{fig: T39window} depicts the window before a cyclic shift, bounded by the vertical grey boundaries, and after a cyclic shift, which is bounded by the vertical blue boundaries.
	
		Consider the union of the first window with its cyclic shift, as depicted in Figure \ref{fig: T39window}. A framed sheaf restricted to this union is determined by vectors $\{v_1,v_2,v_3,v_4, v_5\}$, which are placed at the regions bounded by the first and second strands. In the diagrams in Figure \ref{fig: T39window} the (stalk of) the sheaf is specified in each open region given by the stratification of the front diagram, by associating the vector space spanned by the vectors written in the region. The volume form in each region is given by the ordered wedge product of vectors in that region.

Now we focus on the grey window. Its boundary underlines the two complete flags $$0\subset \langle v_1 \rangle \subset \langle v_1,v_2 \rangle \subset V,\quad  0\subset \langle v_4 \rangle \subset \langle v_4,v_5 \rangle \subset V.$$
Each flag is shared by two nearby windows. To reduce this replication, one can break the symmetry by choosing one flag for each window. Without loss of generality, we choose the flag on the left boundary of each window. In particular, the sheaf restricted to the grey window is reduced to the data of three vectors $\{v_1,v_2,v_3\}$ in $V$.

Note that even though the subspaces $\langle v_4\rangle$ and $\langle v_4,v_5\rangle$ cannot be computed from $\{v_1,v_2,v_3\}$, they are uniquely determined by the next window, and the sheaf is still well-defined over the grey window. We now perform the descending Reidemeister move III depicted in the middle of Figure \ref{fig: T39window}. This R3$^d$ move creates a region and a new vector $u_1$, as we described in the discussion preceding Proposition \ref{prop:EffectSigma1}. From the front, the microlocal support condition for our constructible sheaf implies that
$$\langle v_1,v_2 \rangle = \langle v_2,u_1 \rangle,\qquad \langle v_3,v_4 \rangle = \langle u_1,v_2 \rangle.$$
Hence $u_1$ lies in both $\langle v_1,v_2 \rangle $ and $\langle v_3,v_4 \rangle $. Moreover, the crossing condition at the crossing depicted in red in Figure \ref{fig: T39window} yields that the complex
$$0 \longrightarrow \langle u_1 \rangle \longrightarrow \langle v_1,v_2 \rangle \oplus \langle v_3,v_4 \rangle \longrightarrow V \rightarrow 0$$
is a short exact sequence of $\C$-vector spaces. Therefore
$$\langle u_1 \rangle = \langle v_1,v_2 \rangle \cap \langle v_3,v_4 \rangle,$$
and $u_1$ is the unique vector such that
$$v_1v_2 = v_2u_1.$$ 

This establishes the description of $u_1$ in the statement of Proposition \ref{fig: T39window}. Let us now shift to the blue window. The constructible sheaf restricted to this window is determined by $\{v_2,u_1,v_3,v_4\}$. The fourth vector $v_4$ disappears upon performing the ascending Reidemeister III move, as depicted in the bottom of Figure \ref{fig: T39window}. After this R3$^a$ move, the sheaf is uniquely determined by $\{v_2,u_1,v_3\}$. The morphism induced by $\Sigma_1$ thus starts with
$$(v_1,v_2,v_3) \mapsto (v_2,u_1,v_3),$$
where $u_1 \in \langle v_1,v_2\rangle\cap \langle v_3,v_4\rangle$ and $v_1v_2 = v_2u_1$, and continues to remove $v_4$. These two moves are preceded by the cyclic shift, and their composition yields the expression in the first, and thus any, window in the statement of Proposition \ref{prop:EffectSigma1}, as required.\hfill$\Box$

\subsection{Comments on the Proof} This concludes the proof of Theorem \ref{thm:main}. Before proceeding with Theorem \ref{thm:main2}, the following comments might be clarifying. The geometric loops $\Sigma_1,\delta^2$ are studied in the above proof of Theorem \ref{thm:main} by analyzing their action on the ring of functions $\C[\SM(\La(\beta),\tau)]$ of the framed moduli space $\SM(\La(\beta),\tau)$. It should be equally possible to deduce Theorem \ref{thm:main} by studying their monodromy invariants in the ring of regular functions $\C[\SM(\La(\beta))]$ of the moduli spaces of sheaves, with no frame $\tau$ chosen, with corresponding $(\C^*)^9$-equivariant condition added. Indeed, the positroid embedding of $\SM(\La(\beta),\tau)$ inside the Grassmannian $\Gr(3,9)$ yields an embedding of the moduli of sheaves $\SM(\La(\beta))$ into the quotient $\Gr(3,9)/(\C^*)^9$ of the Grassmannian $\Gr(3,9)$ by the diagonal subgroup of $\GL_9(\C)$ acting on the right, i.e. by column $\C^*$-rescaling.

It is our aesthetic opinion that working directly in the unquotiented Grassmannian $\Gr(3,9)$ yields a clearer understanding of the geometry, thus our choice of using the moduli space of framed sheaves. In terms of cluster algebras, the quotient $\Gr(3,9)/(\C^*)^9$ has no frozen cluster variables, whereas the Grassmannian $\Gr(3,9)$ \cite{FominZelevinsky_ClusterI,ScottGrassmannian} has the cyclically consecutive Pl\"ucker coordinates as frozen cluster variables.

\begin{remark}
The articles \cite{Fraser,ShenWeng} respectively use the affine cone on the projective Grassmannian $\Gr(3,9)$ \cite[Section 3]{Fraser} and the decorated Grassmannian $\mathscr{G}r(3,9)$ \cite[Section 2.1]{ShenWeng}. These can be equivalently considered \cite[Lemma 2.6]{ShenWeng} and correspond to matrices $\mbox{Mat}_{3,9}$ up to the left action of $\SL_3(\C)$, rather than $\GL_3(\C)$, which would yield the projective Grassmannian $\Gr(3,9)$. In terms of the moduli space of framed sheaves $\SM(\La(\beta),\tau)$ used in our proof of Theorem \ref{thm:main}, we should require the additional data of a trivialization of the microlocal monodromy along $\La(\beta)$ itself \cite[Section 5.1]{STZ}. By context, it seems appropriate to refer to this space as the moduli space of decorated sheaves. The line of argument above should also work by using the decorated positroid embedding of the space of decorated sheaves into the decorated Grassmannian.\hfill$\Box$
\end{remark}

Let us now move forward with Theorem \ref{thm:main2}. Note that Theorem \ref{thm:main} on its own allows us to conclude Corollaries \ref{cor:infinitefillings} and \ref{cor:weinstein} in the cases $(n,m)\in\SH\setminus\{(4,4),(4,5),(5,5)\}$, and Corollaries \ref{cor:concordances} and \ref{cor:fundamentalgroup} for $\La(3,6)$. In order to cover the Legendrian links $\La(4,4)$, $\La(4,5)$ and $\La(5,5)$, and for completeness, we now include the proof of Theorem \ref{thm:main2}, which is in line with that of Theorem \ref{thm:main} above.

\section{The $M_{0,4}$ representation for $\La(4,4)$}\label{sec:mainproof2}

In this section we prove Theorem \ref{thm:main2}. The argument reproduces the strategy for Theorem \ref{thm:main} above. In this case, the braid is $\beta=(\sigma_1\sigma_2\sigma_3)^8$ and the moduli space $\SM(\La(\beta))$ is identified with a positroid cell $\Pi_{r(\beta)}\sse\Gr(4,8)$ by the same procedure. The action of the Legendrian loops $\Xi_1,\Xi_2,\Xi_3$ is described by the following three crucial Propositions:

\begin{prop}\label{prop:EffectXi1} The Legendrian loop $\Xi_1$ induces the morphism
		$$(v_1,v_2,v_3,v_4;v_5,v_6,v_7,v_8)\longmapsto(v_2,u_1,v_3,v_4;v_6,u_2,v_7,v_8),$$
		where $u_1,u_2$ are given by the intersections
		$$u_1\in\langle v_1,v_2\rangle\cap \langle v_3,v_4,v_5\rangle,\quad u_2\in\langle v_5,v_6\rangle\cap \langle v_7,v_8,v_1\rangle,$$
		and the normalizing conditions $v_1v_2=v_2u_1$ and $v_5v_6=v_6u_2$ in $V\wedge V$.
\end{prop}

\begin{prop}\label{prop:EffectXi2} The Legendrian loop $\Xi_2$ induces the morphism
	$$(v_1,v_2,v_3,v_4;v_5,v_6,v_7,v_8)\longmapsto(v_1,v_3,u_1,v_4;v_5,v_7,u_2,v_8),$$
	where $u_1,u_2$ are given by the intersections
	$$u_1\in\langle v_2,v_3\rangle\cap \langle v_4,v_5,v_6\rangle,\quad u_2\in\langle v_6,v_7\rangle\cap \langle v_8,v_1,v_2\rangle,$$
	and the normalizing conditions $v_2v_3=v_3u_1$ and $v_6v_7=v_7u_2$ in $V\wedge V$.
\end{prop}

\begin{prop}\label{prop:EffectXi3} The Legendrian loop $\Xi_3$ induces a morphism
	$$(v_1,v_2,v_3,v_4;v_5,v_6,v_7,v_8)\longmapsto(v_1,v_2,v_4,u_1;v_5,v_6,v_8,u_2),$$
	where $u_1,u_2$ are given by the intersections
	$$u_1\in\langle v_3,v_4\rangle\cap \langle v_5,v_6,v_7\rangle,\quad u_2\in\langle v_7,v_8\rangle\cap \langle v_1,v_2,v_3\rangle,$$
	and the normalizing conditions $v_3v_4=v_4u_1$ and $v_7v_8=v_8u_2$ in $V\wedge V$.
\end{prop}

Propositions \ref{prop:EffectXi1}, \ref{prop:EffectXi2} and \ref{prop:EffectXi3} are proven at the end of this section. The action of the group $\Gamma_2=\langle\Xi_1,\Xi_2,\Xi_3\rangle$ in the set of 8-tuples of vectors, representing a point in $\Gr(4,8)$, yields via pull-back an action on a subset of the homogeneous coordinate ring $\C[\Gr(4,8)]$. For the braid $\La(4,4)$ it does not suffice to study the $\Gamma_2$-orbit of a Pl\"ucker coordinate, as we directly did for Theorem \ref{thm:main}, but rather a set of Pl\"ucker coordinates. In this proof for Theorem \ref{thm:main2}, we directly refer to known algebraic arguments whose nature is on par with Subsection \ref{ssec:FaithfulPSL}, as follows. Indeed, \cite[Lemma 10.8]{Fraser} proves that the group $\langle \Xi_1,\Xi_2,\Xi_3\rangle$ generated by the monodromies of the three Legendrian loops generates a faithful action of $\PSL(2,\Z)\cong\Z_2*\Z_3$ on the (cluster) automorphism group of the coordinate ring $\C[\Gr(4,8)]$. This is achieved by studying the orbit of the Pl\"ucker set:
$$\mathscr{P}=\{P_{1378},P_{2348},P_{2367},P_{4678},P_{3457},P_{2347},P_{2378},P_{3678},P_{3467}\},$$
which is a cluster seed for a triangulation of the annulus with four boundary marked points. By \cite[Proposition 2.7]{FarbMargalit12}, the mapping class group of the four-punctured sphere $M_{0,4}$ is isomorphic to the semidirect product $\PSL(2,\Z)\ltimes(\Z_2\times\Z_2)$. The article \cite[Theorem 9.14]{Fraser} also shows that this faithful action of $\PSL(2,\Z)$ extends to the mapping class group $M_{0,4}$ as required. This is achieved explicitly by studying the four cosets of $\PSL(2,\Z)$ into $M_{0,4}$. In the algebraic argument the set $F$ can be chosen to be the union of four sets, as follows. The first set $S$ is the union of a finite number of cluster charts \cite[Section 10.2]{Fraser} containing the set $\mathscr{P}$ of Pl\"ucker coordinates above, and the remaining three sets are the coset translates $\Xi_3(S)$, $\Xi_3\Xi_2(S)$, and $\Xi_3\Xi_2\Xi_1(S)$. Here	 $\Xi_3,\Xi_3\Xi_2$ and $\Xi_3\Xi_2\Xi_1$ are each a right coset representative for each of the three non-trivial cosets of the inclusion of $\PSL(2,\Z)$ into $M_{0,4}$ above.\hfill$\Box$

The crucial ingredient for the proof of Theorem \ref{thm:main2} above is the statement that the Legendrian loops we constructed in Section \ref{sec:loops} indeed induce an action of the (spherical) braid group $B_4$. This is precisely the content of Propositions \ref{prop:EffectXi1}, \ref{prop:EffectXi2} and \ref{prop:EffectXi3}, which describe the algebraic effect of the Legendrian loops $\Xi_1,\Xi_2$ and $\Xi_3$. Let us now prove these three propositions.

\subsection{Proof of Proposition \ref{prop:EffectXi1}} Let us consider the braid words $\beta = (\sigma_1\sigma_2\sigma_3)^8$ and $\beta_0 = (\sigma_1\sigma_2\sigma_3)^4$. Following the notation in the proof of Proposition \ref{prop:EffectSigma1} above, each $\beta_0$ is a window and $\beta = \beta_0^2$ is the concatenation of two windows. Similar to the proof of Proposition \ref{prop:EffectSigma1}, it suffices to compute the induced morphism in a window.

\begin{center}
	\begin{figure}[h!]
		\centering
		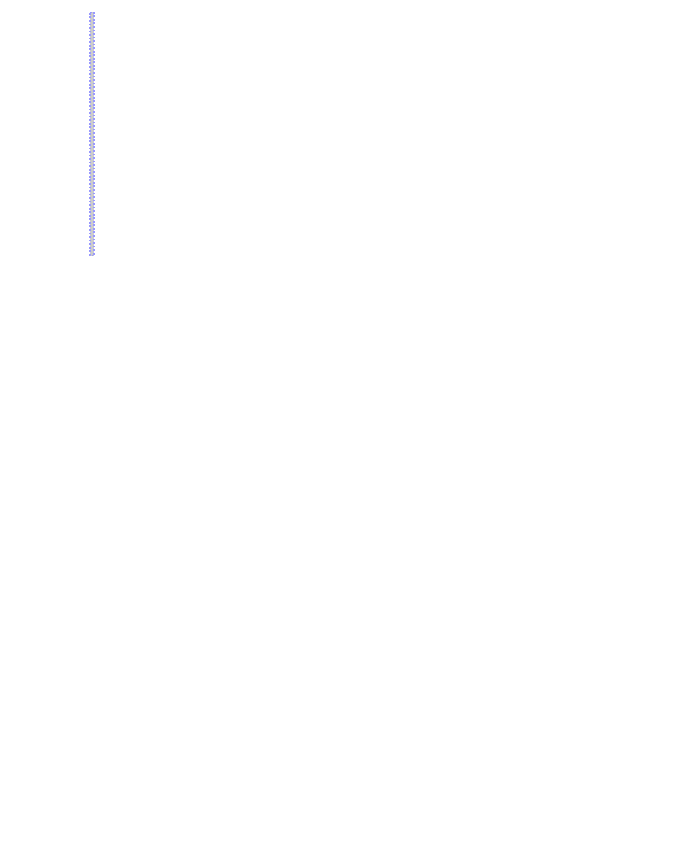
		\caption[]{First window $\beta_0$ of the Legendrian braid $\beta=(\sigma_1\sigma_2\sigma_3)^8$ associated to the Legendrian loop $\Xi_1$ for the Legendrian link $\Lambda(4,4)$. There are a total of two windows.}
		\label{fig:T48window}
	\end{figure}
\end{center}

Consider the union of the first window and its one-term cyclic shift, depicted in the top diagram of Figure \ref{fig:T48window}. The window before the shift has grey boundaries. We choose to include the sheaf data on the left boundary in this window, and leave the sheaf data on the right boundary to the next window. With this choice, a framed constructible sheaf in the grey window is determined by four vectors $\{v_1,v_2,v_3,v_4\}$ in $V$.

Now we study the morphism induced by the Legendrian loop $\Xi_1$. The sequence of braid moves can be carried out as the concatenation of two Legendrian isotopies $\Psi_1^{(2)}\circ \Psi_1^{(1)}$. These two Legendrian isotopies $\Psi_t^{(1)}$ and $\Psi_t^{(2)}$, $t\in[0,1]$, are defined in Section \ref{sec:loops}. In Figure \ref{fig:T48window}, $\Psi^{(1)}_t$ corresponds to the Legendrian isotopy from the top diagram to the middle diagram, and $\Psi^{(2)}_t$ is depicted from the middle diagram to the bottom diagram.
After performing the Legendrian isotopy $\Psi_t^{(1)}$, $t\in[0,1]$, the diagram introduces a new vector $u_1$. From the diagram, we see that $\langle v_1, v_2 \rangle = \langle v_2, u_1 \rangle$ and $\langle u_1,v_3, v_4\rangle = \langle v_3,v_4,v_5 \rangle$. Hence $u_1\in \langle v_1, v_2 \rangle  \cap \langle v_3,v_4,v_5 \rangle$. To argue that the intersection is a one-dimensional subspace, we should discard the case that $\langle v_1, v_2 \rangle  \subset \langle v_3,v_4,v_5 \rangle$. Inside the middle figure, the condition at the red crossing yields a short exact sequence of complex vector spaces:
$$0\rightarrow \langle u_1, v_3 \rangle \rightarrow \langle v_1, v_2, v_3 \rangle \oplus \langle v_3, v_4, v_5 \rangle \rightarrow V\rightarrow 0.$$
If $\langle v_1, v_2 \rangle$ is contained in $\langle v_3,v_4,v_5 \rangle$, so is $\langle v_1, v_2, v_3 \rangle$. Then it is impossible to map the direct sum onto $V$, which is a contradiction. Therefore
$$\langle u_1 \rangle= \langle v_1, v_2 \rangle  \cap \langle v_3,v_4,v_5 \rangle,$$
and the vector $u_1$ can be uniquely determined by 
$$v_1v_2 = v_2u_1.$$

At this stage, there are five vectors $\{v_2,u_1,v_3,v_4,v_5\}$ inside the (blue) shifted window. There is a redundancy which is removed via the Legendrian isotopy $\Psi^{(2)}_t$. The bottom diagram in Figure \ref{fig:T48window} specifies how to determine the constructible sheaf using the four vectors $\{v_2,u_1,v_3,v_4\}$. In the end, the only regions including $v_5$ are connected to the right blue boundary, which is determined by the next window. Iterating this procedure in each window, we obtain that the morphism determined by $\Xi_1$ is indeed that of the statement of Proposition \ref{prop:EffectXi1}.\hfill$\Box$

\subsection{Proof of Proposition \ref{prop:EffectXi2}}
\begin{center}
	\begin{figure}[h]
		\centering
		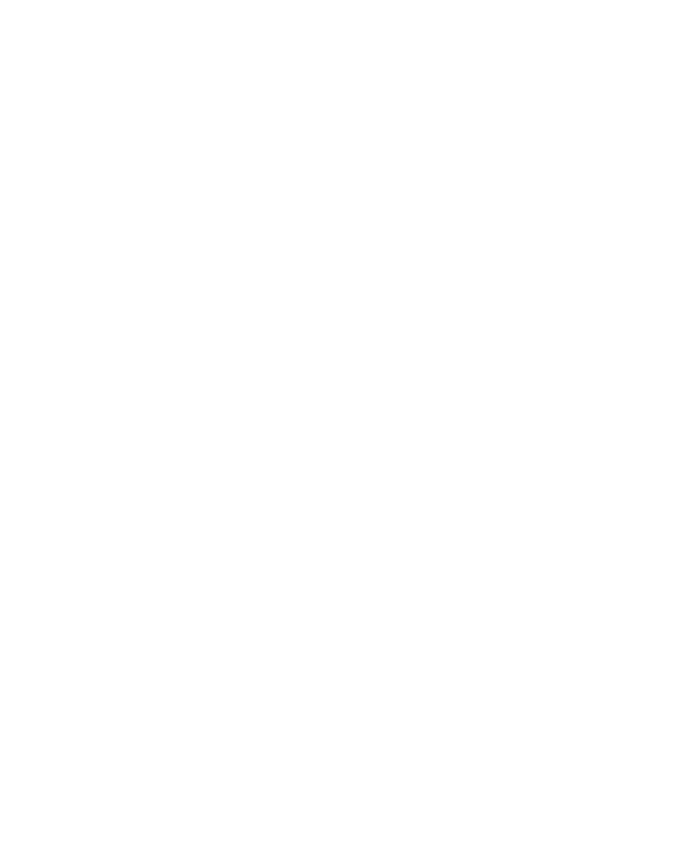
		\caption[]{Window for the Legendrian loop $\Xi_2$ on $\Lambda(4,4)$ with the needed information for the proof of Proposition \ref{prop:EffectXi2}.}
		\label{fig:Xi2T48window}
	\end{figure}
\end{center}

Let us consider the Legendrian isotopy $\Psi^{(1)}_t$, as defined in Section \ref{sec:loops}. This is the first of two pieces which constitute the Legendrian loop $\Xi_2$. This Legendrian isotopy $\Psi^{(1)}_t$ is depicted from the top to the middle in Figure \ref{fig:Xi2T48window}; we have labeled two of the regions in the middle picture, each being assigned a $2$-dimensional vector space, denoted $W_1$ and $W_2$. Note that $W_1 = \langle v_2, v_3\rangle$, since this region already exists in the front at the top row of Figure \ref{fig:Xi2T48window}. The second vector space $W_2$ is given by the intersection $\langle v_2,v_3,v_4 \rangle \cap \langle v_4,v_5,v_6 \rangle $, following the condition at the red crossing in the middle picture. This determines the algebraic effect of the Legendrian isotopy $\Psi^{(1)}_t$.

Let us continue with the second Legendrian isotopy $\Psi^{(2)}_t$. This Legendrian isotopy creates a new region with a vector $u_1$, as depicted in Figure \ref{fig:Xi2T48window}. The vector spaces $W_1$ and $W_2$ can then be described by using the vector $u_1$. Indeed, we have $W_1 = \langle v_3,u_1\rangle$ and $W_2 =  \langle u_1,v_4\rangle $. An argument in line with that of the proof of Proposition \ref{prop:EffectXi1} concludes that $u_1\in\langle v_2,v_3\rangle\cap \langle v_4,v_5,v_6\rangle$ and it is uniquely determined by $v_2v_3=v_3u_1$, as required. This concludes the desired transformation for the first window. The transformations for the remaining windows are concluded similarly.\hfill$\Box$

\begin{center}
	\begin{figure}[h!]
		\centering
		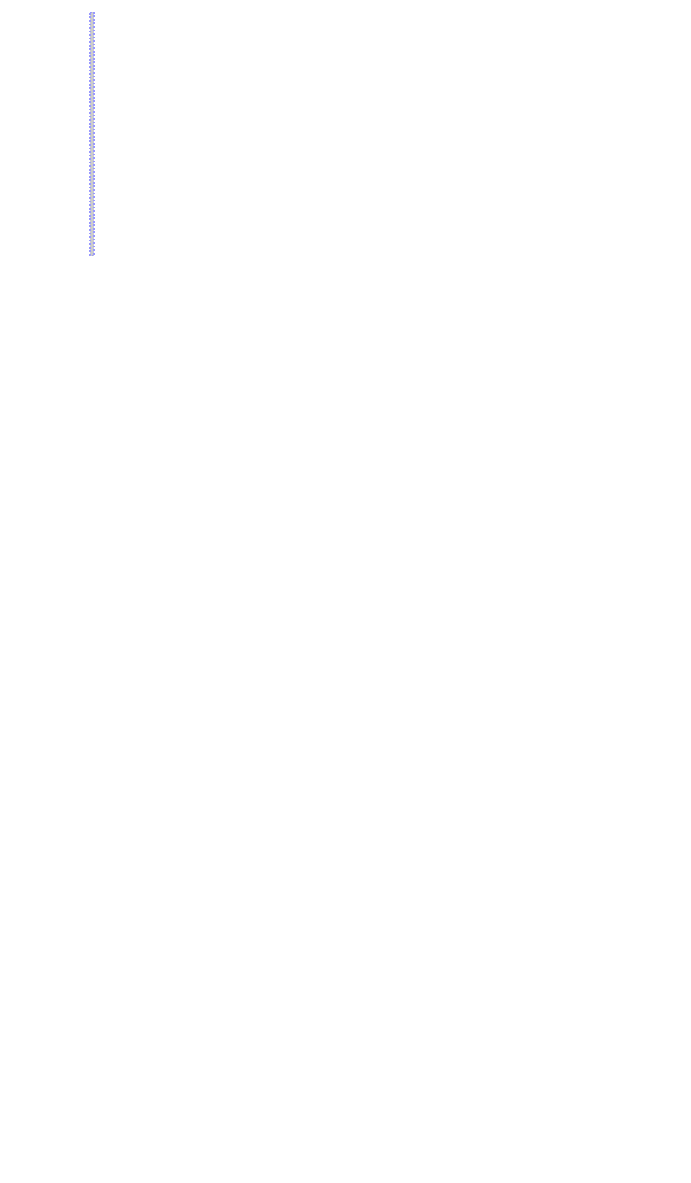
		\caption[]{Windows for the Legendrian loop $\Xi_3$ on $\Lambda(4,4)$ containing the required information for the proof of Proposition \ref{prop:EffectXi3}.}
		\label{fig:Xi3T48window}
	\end{figure}
\end{center}

\subsection{Proof of Proposition \ref{prop:EffectXi3}} The argument is identical to that in Propositions \ref{prop:EffectXi1} and \ref{prop:EffectXi2}, and thus we only provide the core steps. In particular, we have depicted the Legendrian loop $\Xi_3$ in Figure \ref{fig:Xi3T48window} as well as its effect in three different pieces $\Xi_3 = \Psi_1^{(3)} \circ \Psi_1^{(2)} \circ \Psi_1^{(1)} $, as recorded in Section \ref{sec:loops}. In short, the core information in studying the effect of $\Xi_3$ can be described as follows:

\begin{itemize}
	\item[-] After $\Psi_1^{(t)}$, the subspaces are uniquely determined as indicated in the figure. The vector space spanned by $v_2$ disappears but the vector $v_2$ can be recovered from the new data. Namely, it is determined by the intersection of $\langle v_1, v_2\rangle$ and $\langle v_2,v_3,v_4\rangle$, both of which are stalks of some regions in the front diagram, and the volume form in either one of these vector spaces.\\
	
	\item[-] After $\Psi_2^{(t)}$, the subspaces are also uniquely determined as indicated. The data of $v_2$ remains in the diagram implicitly.\\
	
	\item[-] The Legendrian isotopy $\Psi_3^{(t)}$, pulls down two strands in Figure \ref{fig:Xi3T48window}  which are colored in red. The red strand on the left recovers the vector $v_2$. The red strand on the right introduces a new vector $u_1$, which satisfies $v_3v_4 = v_4u_1$, and $u_1v_5v_6 = v_5v_6v_7$. By a similar argument with that for $\Xi_1$ and $\Xi_2$, we see that $u_1 \in \langle v_3,v_4\rangle\cap \langle v_5,v_6,v_7\rangle$ and that $v_3v_4 = v_4u_1$. 
	
\end{itemize}
In conclusion, the morphism sends the first window from the 4-tuple $(v_1,v_2,v_3,v_4)$ to the 4-tuple $(v_1,v_2,v_4,u_1)$ as required. The second window is concluded in similar manner.\hfill$\Box$

\section{Corollaries and Applications}\label{sec:corollaries}

In this section we prove Corollaries \ref{cor:infinitefillings},\ref{cor:lisinfinite}, \ref{cor:concordances}, \ref{cor:fundamentalgroup} and \ref{cor:weinstein}.

First, Corollary \ref{cor:concordances} follows by observing that a trivial concordance in the Lagrangian concordance monoid $\mathbb{L}(3,6)$, and $\mathbb{L}(4,4)$, would induce a trivial map on $\cM(\La(3,6))$, and $\cM(\La(4,4))$ respectively. Theorems \ref{thm:main} and \ref{thm:main2} imply that the loops $\Sigma_1,\delta^2$, for $\La(3,6)$ and $\Xi_1,\Xi_2,\Xi_3$, for $\La(4,4)$, induce Legendrian loops which act non-trivially on $\cM(\La(3,6))$, and $\cM(\La(4,4))$ respectively. Hence the concordances induced by graphing these Legendrian loops are themselves non-trivial. The same argument concludes Corollary \ref{cor:fundamentalgroup}.

\begin{center}
	\begin{figure}[h!]
		\centering
		\includegraphics[scale=0.7]{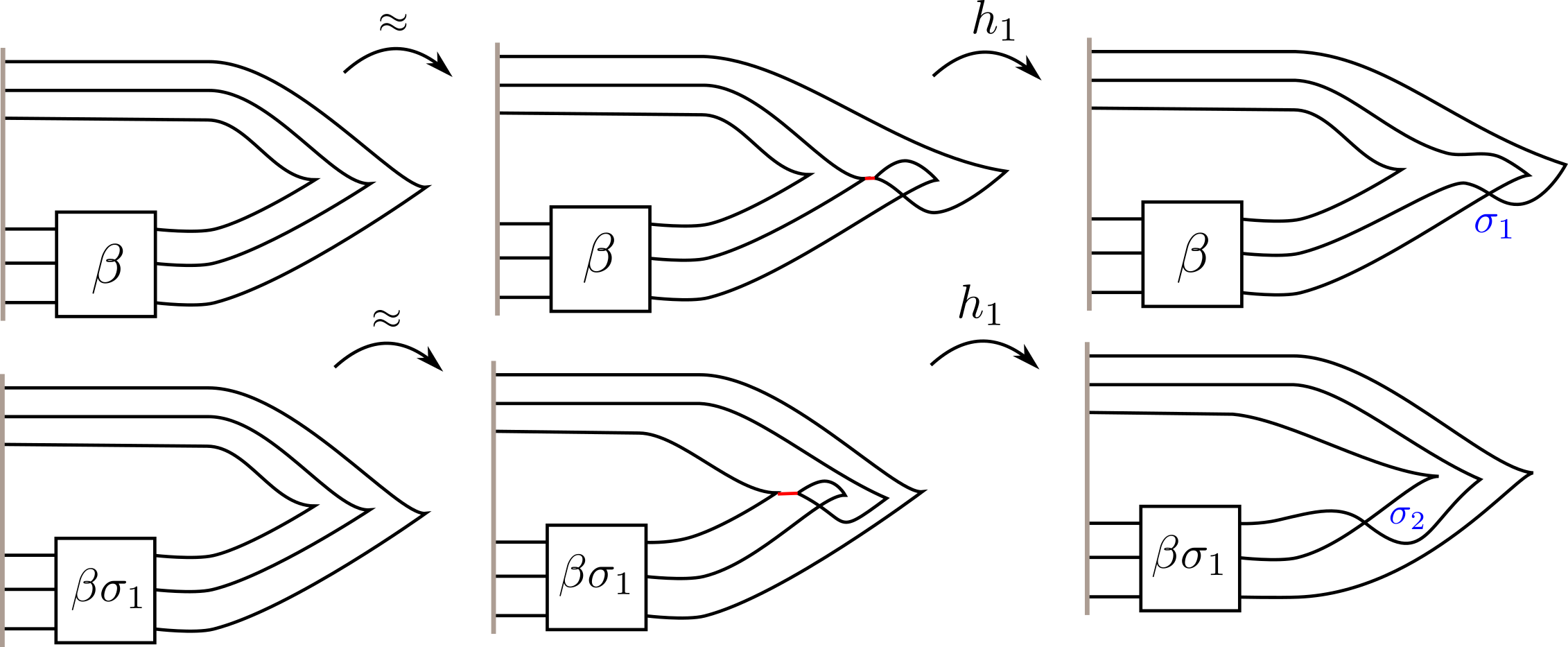}
		\caption{Exact Lagrangian Cobordism from $\beta$ to $\beta\sigma_1\sigma_2$.}
		\label{fig:LegLinks_Cobordism}
	\end{figure}
\end{center}

Let us now address Corollary \ref{cor:infinitefillings} and Corollary \ref{cor:weinstein}, which shall follow from Theorems \ref{thm:main} and \ref{thm:main2}, with the addition of the upcoming Proposition \ref{prop:cobordism}. For that, let us consider the two-sided closure, i.e.~the rainbow closure, of the braid word $\beta=(\sigma_1\cdot\sigma_2\cdot\ldots\cdot\sigma_{n-1})^{m}$ as depicted in the upper leftmost diagram in Figure \ref{fig:LegLinks_Cobordism}. Let us denote the Legendrian associated to this front $\La(\beta)$. Corollary \ref{cor:infinitefillings} is proven with the following geometric construction:

\begin{prop}\label{prop:cobordism} Let $\La(n,m)=\La(\beta)$ be the Legendrian torus link given by the braid $$\beta=(\sigma_1\cdot\ldots\cdot\sigma_{n-1})^m.$$
There exists a decomposable Lagrangian cobordism from $\La(n,m)$ to $\La(n,m+1)$ whose Lagrangian handles have isotropic spheres away from the region with the $\beta$-braiding. Similarly, there exists a decomposable Lagrangian cobordism from $\La(n,m)$ to $\La(n+1,m)$.
\end{prop}
	
\begin{center}
	\begin{figure}[h!]
		\centering
		\includegraphics[scale=0.8]{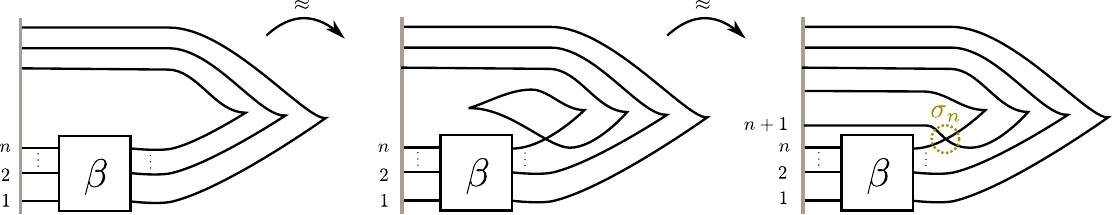}
		\caption{Legendrian isotopy from $\beta\in\mbox{Br}_n$ to $\beta\sigma_n\in\mbox{Br}_{n+1}$.}
		\label{fig:LegLinks_Cobordism2}
	\end{figure}
\end{center}

\begin{proof}
For any given $\sigma_i$, $1\leq i\leq n-1$, it suffices to construct a decomposable Lagrangian cobordism with concave end $\La(\beta)$ and convex end $\La(\beta\sigma_i)$. For that, we first perform an upwards Reidemeister I move on the right lower strand for the $i$th rightmost cusp. Then, the left cusp created in this Reidemeister I move can be isotoped, without introducing crossings in the front, to the same level as the rightmost cusp for the $(i+1)$th strand. This is depicted in the second and fifth diagrams of Figure \ref{fig:LegLinks_Cobordism} in the cases of $\sigma_1$ and $\sigma_2$. Once these two cusps are aligned, we perform a reverse pinched-move \cite{BourgeoisSabloffTraynor15,PanCatalanFillings} allowing this pair of opposite cusps to become two parallel strands. This corresponds to a Lagrangian $1$-handle attachment $h_1$, and it is depicted in the second to third, and fifth to sixth diagrams in Figure \ref{fig:LegLinks_Cobordism}. The decomposable Lagrangian cobordisms just described can be independently and repeatedly performed for different $\sigma_i$, $1\leq i\leq n-1$. In particular, by applying this cobordism for the Artin generators $\sigma_1,\sigma_2$ through $\sigma_{n-1}$, we obtain a decomposable Lagrangian cobordism from $\La(\beta)$ to $\La(\beta(\sigma_1\sigma_2\cdot\ldots\cdot\sigma_{n-1}))$, which implies the statement in the proposition when applied to the braid $\beta=(\sigma_1\cdot\ldots\cdot\sigma_{n-1})^m$.

The decomposable Lagrangian cobordism from $\La(n,m)$ to $\La(n+1,m)$ is built similarly. First, the Legendrian link whose front is the rainbow closure of a $k$-stranded positive braid word $\beta\in\mbox{Br}_k$ is Legendrian isotopic to the Legendrian link whose front is the rainbow closure of a $(k+1)$-stranded positive braid word $\beta\sigma_k\in\mbox{Br}_{k+1}$. This is proven by performing a Legendrian Reidemeister I move, which introduces the $\sigma_k$-crossing, and it is depicted in Figure \ref{fig:LegLinks_Cobordism2}. Thus the front given by the rainbow closure of $(\sigma_1\cdot\ldots\cdot\sigma_{n-1})^m$ is front homotopic to the rainbow closure of $(\sigma_1\cdot\ldots\cdot\sigma_{n-1})^{m-1}\cdot(\sigma_1\cdot\ldots\cdot\sigma_{n-1}\sigma_n)$; these both give the Legendrian $\La(n,m)$, the latter front using a $(n+1)$-stranded braid. Second, it now suffices to add $(m-1)$ new positive crossings $\sigma_n$ to $(\sigma_1\cdot\ldots\cdot\sigma_{n-1}\sigma_n)^{m}$, which again can each be inserted via an index-1 decomposable (exact) Lagrangian cobordism. Note that it is possible to insert any positive crossing in the middle of a braid word $\beta$ (not just at its rightmost end) with such a Lagrangian cobordism. Indeed, one may apply a cyclic shift $\delta^k$, for some $k\in\N$, so that the location where the new crossing is to be inserted is to the right of $\delta^k(\beta)$, then apply the exact Lagrangian cobordism from Figure \ref{fig:LegLinks_Cobordism}, and compose with the inverse of the Legendrian isotopy $\delta^k$. Inserting these $(m-1)$ positive crossings $\sigma_n$ allows us to arrive to $(\sigma_1\cdot\ldots\cdot\sigma_{n-1}\sigma_n)^{m}$ from $(\sigma_1\cdot\ldots\cdot\sigma_{n-1})^{m-1}\cdot(\sigma_1\cdot\ldots\cdot\sigma_{n-1}\sigma_n)$.
This yields the required decomposable Lagrangian cobordism from $\La(n,m)$ to $\La(n+1,m)$.
\end{proof}

Note that Proposition \ref{prop:cobordism} holds for any pair $(n,m)\in\N\times\N$, with no constraint $n\leq m$ nor $(n,m)\in\SH$.

\subsection{Proof of Corollary \ref{cor:infinitefillings}}
Let us first prove that $\La(3,6)$ has infinitely many Lagrangian fillings. Fix an exact Lagrangian filling $L\sse(\R^4,\omega)$ for $\La(3,6)$ obtained via a pinching sequence from the front diagram on the left of Figure \ref{fig:BraidsIntro}. Smoothly, this must be a thrice punctured genus-3 surface. By \cite[Proposition 2.15]{STWZ}, this exact Lagrangian filling yields an open inclusion $\mbox{Loc}(L)\sse\SM(\La(3,6))$, where Loc$(L)$ denotes the space of framed local systems in $L$. Now, given a Legendrian loop $\vartheta\in\Gamma=\langle[A],[B]\rangle$ in the group generated by the Legendrian loops $A,B$ (or equivalently $\Sigma_1,\delta^2$), we consider the Lagrangian filling $L_\vartheta$ obtained by applying the Legendrian loop $\vartheta$ to the Legendrian $\La(3,6)$ and then performing the fixed pinching sequence for the Lagrangian filling fixed above. Choose an infinite sequence of distinct elements $(\vartheta_i)_{i\in\N}\in\PSL(2,\Z)$, since $\vartheta_i$ are distinguished by their action on the infinite cluster charts of $\Gr(3,9)$, the inclusions $\mbox{Loc}(L_{\vartheta_i})\sse\SM(\La(3,6))$ yield infinitely many distinct cluster charts. In consequence, the Lagrangian fillings $L_{\vartheta_i}$ are not Hamiltonian isotopic \cite[Proposition 6.1]{STWZ}. The same argument holds for the Legendrian link $\La(4,4)$ once we use the representation in Theorem \ref{thm:main2} and the mapping class group $M_{0,4}$.

Let $\La(n,m)$ be given with $(n,m)\in\SH$ different from $(4,4),(4,5),(5,5)$. The construction in the proof of Proposition \ref{prop:cobordism} yields a decomposable Lagrangian cobordism from $\La(3,6)$ to $\La(n,m)$. This exact Lagrangian cobordism yields an injective map between the equivalence classes of objects of the associated $\Aug_+$ categories, i.e. distinct augmentations up to isomorphism (including DGA homotopy) for $\La(3,6)$ yield, upon composing with the DGA map induced by this Lagrangian cobordism, distinct augmentations for $\La(n,m)$. Injectivity in the case of knots is proven in \cite[Theorem 1.5]{PanCatalan_AugmentationMap}, the case of links is analogous and it is detailed in \cite{CLLMPT}; see also Remark \ref{rmk:injective} below. Since there are infinitely many Lagrangian fillings for $\La(3,6)$ distinguished by their sheaves, the correspondence between augmentations and sheaves \cite[Theorem 1.3]{NRSSZ} implies that these Lagrangian fillings are distinguished by their augmentations on their Chekanov-Eliashberg algebra \cite{Chekanov02}. Thus, the infinitely many Lagrangian fillings of $\La(3,6)$ concatenated with the Lagrangian cobordism in Proposition \ref{prop:cobordism} induce non-isomorphic augmentations for $\La(n,m)$. In consequence, the infinitely many Lagrangian fillings of $\La(3,6)$ yield infinitely many Lagrangian fillings of $\La(n,m)$. For the remaining case of $(n,m)=(4,5),(5,5)$ we apply Proposition \ref{prop:cobordism} to obtain a cobordism from $\La(4,4)$ to $\La(4,5)$ or $\La(5,5)$ and proceed identically.\hfill$\Box$

\begin{remark}\label{rmk:injective}
Let $\La_-,\La_+$ be two Legendrian links such that there exists a decomposable exact Lagrangian cobordism from $\La_-$ to $\La_+$. In these hypotheses, the argument for Corollary \ref{cor:infinitefillings} above uses the following fact: if $\La_-$ is a Legendrian link that admits infinitely many Lagrangian fillings which are distinguished by augmentations (resp.~by sheaves)\footnote{E.g. they induce different sheaves in the analogous category $\mathcal{C}_1(\La_-)$, in the notation of \cite{NRSSZ}.} -- e.g. they yield different objects in the $\Aug_+$ category -- then $\La_+$ is a Legendrian link that also admits infinitely many Lagrangian fillings which are distinguished by augmentations (resp.~by sheaves).

As mentioned, in the case of augmentations and both $\La_-,\La_+$ knots, this fact is known to hold for an arbitrary exact Lagrangian cobordism, not necessarily decomposable, by \cite[Theorem 1.5]{PanCatalan_AugmentationMap}. Nevertheless, a much simpler argument exists if one assumes that the exact Lagrangian cobordism is decomposable, as it is in our case. Then \cite[Proposition 7.5]{CasalsNg} shows that this fact is true, now also including the general case where both $\La_-,\La_+$ are allowed to be links, which suffices for our purposes.\hfill$\Box$
\end{remark}

The cluster modular groups of the remaining Grassmannians $\Gr(n,m+n)$, with the pair $(n,m)\in(\N\times\N)\setminus\SH$, are known to be finite \cite{ASS,Fraser}. Thus, for these remaining Legendrian links $\La(n,m)$, $(n,m)\in(\N\times\N)\setminus\SH$, our arguments will only yield a representation of a finite group. In particular, we are almost certain that our results are sharp, i.e.~we conjecture that the Legendrian torus links $\La(n,m)$ have finitely many Lagrangian fillings if $(n,m)\not\in\SH$. In fact, we believe that the Legendrian torus links $\La(2,n)$ must have exactly $\frac{1}{n+1}{2n\choose n}$ Lagrangian fillings, $\La(3,3)$ should have exactly $50$ Lagrangian fillings, and $\La(3,4)$ and $\La(3,5)$ will have exactly 883 and 25080 Lagrangian fillings respectively.
	
\begin{remark}
The numbers 50, 883 and 25080 are the number of cluster seeds for the finite type cluster algebras of types $D_4$, $E_6$ and $E_8$, respectively. See \cite[Proposition 3.8]{FominZelevinsky03}, \cite[Theorem 1.13]{FominZelevinsky_ClusterII}, and \cite[Section 5]{CasalsADE}. Note that these numbers are strictly greater than the number of corresponding maximal pairwise weakly separated collections, and thus each correspondingly greater than the number of embedded exact Lagrangian fillings constructed in \cite[Proposition 6.2]{STWZ}. For instance, \cite{STWZ} builds 34 exact Lagrangian fillings for $D_4$ (resp.~259 for $E_6$), namely those corresponding to maximal pairwise weakly separated collections with $k=3$ and $n+k=6$ (resp.~ $k=3$ and $n+k=7$). Yet, the remaining 16 (resp.~574) clusters of $\Gr(3,6)$ (resp.~$\Gr(3,7)$), are also inhabited by embedded exact Lagrangian fillings; see \cite{CasalsADE} and references therein.\hfill$\Box$
\end{remark}

\subsection{Proof of Corollary \ref{cor:lisinfinite}} Let $\La\sse(S^3,\xi_\st)$ be any Legendrian link with an exact Lagrangian cobordism $\La(3,6)\preceq\La$, or $\La(4,4)\preceq\La$. The argument for Corollary \ref{cor:infinitefillings} implies that $\La$ itself will have infinitely many exact Lagrangian fillings. This readily implies Corollary \ref{cor:lisinfinite}. Indeed, by \cite[Theorem 1.1]{BTZ_HyperbolicKnots} the twisted torus knots $K_{p,q,r,s}=T(p,q,r,s)$ with $1<r<p<q$ and $18p<s$ are hyperbolic knots. Let $\La_{K_{3,7,2,s}}$ be the maximal-tb Legendrian representative obtained from the positive braid associated to the $T$-knot $T(3,7,2,s)$, with $54<s$ and $s$ even, as described in \cite[Section 1]{Birman_Lorenz}. Then there exists an exact Lagrangian cobordism $\La(3,6)\preceq\La_{K_{3,7,2,s}}$, and hence $\La_{K_{3,7,2,s}}$ is a hyperbolic knot which admits infinitely many exact Lagrangian fillings. The same argument applies to the twisted torus links $T(p,q,kq,s)$, $p,q,k,s\in\N$, which are proven to be $(q,p+k^2qs)$-cables of the torus knot $T(k,ks+1)$ in \cite{Lee_SatelliteKnots}.\hfill$\Box$

The argument above can be applied in a more ad hoc manner to show that certain knots have Legendrian representatives with infinitely many fillings. For instance, the hyperbolic knot $K=k(4_3)$, which is one of the simplest hyperbolic knots (with four ideal teatrahedra in its complement \cite{CDW_SimplestHyperbolic}) is the twisted torus knot $T(3,8,2,1)$. Given that there exists an exact Lagrangian cobordism $\La(3,8)\preceq\La_{K_{3,8,2,1}}$, and $\La(3,8)$ admits infinitely many exact Lagrangian fillings, we have that the Legendrian knot $\La_{K_{3,8,2,1}}$, which is smoothly $k(4_3)$, also admits infinitely many exact Lagrangian fillings.

\subsection{Proof of Corollary \ref{cor:weinstein}}
Consider an infinite collection of the exact Lagrangian fillings $\{L_i\}_{i\in\N}$ constructed in Corollary \ref{cor:infinitefillings}, and denote by $\overline{L}_i\sse M(n,m)$ the exact Lagrangian surfaces obtained by capping $L_i$ with the unique defining $2$-handle of $M(n,m)$, $i\in\N$. 
By the equivalences between sheaves and augmentations \cite[Theorem 1.3]{NRSSZ}, these Lagrangian fillings $\{L_i\}_{i\in\N}$ are distinguished by the augmentations they induce in the Chekanov-Eliashberg differential graded algebra $A_{n,m}$ of $\La(n,m)$. The wrapped Fukaya categories of the Weinstein manifolds $M(n,m)$ are generated by their respective unique cocore $C$ of their defining $2$-handle \cite{Abouzaid,CRGG}, i.e. the wrapped Fukaya category is identified with the category of dg-modules over $\mbox{End}(C,C)\cong A_{n,m}$. Hence, the Lagrangian surfaces $\overline{L}_i$, whose wrapped Floer complex $WF(C,\overline{L}_i)$ has a unique generator, yield distinct 1-dimensional $A_{n,m}$-modules. Thus $\{\overline{L}_i\}_{i\in\N}$ represent distinct objects in the wrapped Fukaya category and $\{\overline{L}_i\}_{i\in\N}$ are an infinite collection of pairwise non-Hamiltonian isotopic exact Lagrangians.\hfill$\Box$


\bibliographystyle{plain}
\bibliography{LagrFillings}

\end{document}